\documentclass[12pt,twoside]{article}
\usepackage[hmargin=0.8in,vmargin=0.9in]{geometry}
\usepackage{fancyhdr}
\usepackage{graphicx}
\usepackage{subfigure}
\usepackage{amssymb}
\usepackage{amsmath}
\usepackage{amsfonts}
\usepackage{theorem}
\usepackage{mathrsfs}
\usepackage{mathtools}
\usepackage{bm}
\usepackage{color}
\usepackage{setspace}
\usepackage{exscale}
\usepackage{relsize}
\usepackage{epstopdf}
\usepackage{float}
\usepackage{cite}
\usepackage{nicefrac}

\usepackage{booktabs,multirow} 
\usepackage{array} 
\usepackage{paralist} 
\usepackage{verbatim} 
\usepackage{subfigure} 

\theoremstyle{plain}			
\newtheorem{thm}{Theorem}[section]

\newtheorem{rmk}[thm]{Remark}
{\theorembodyfont{\rmfamily}\newtheorem{remark}{Remark}[section]}

\usepackage{tikz}
\usetikzlibrary{positioning}
\usepackage{xcolor}

\allowdisplaybreaks[1]

\numberwithin{equation}{section}
\numberwithin{figure}{section}
\numberwithin{table}{section}

\newcommand\eref[1]{(\ref{#1})}

\newcommand*\xbar[1]{%
  \hbox{%
    \vbox{%
      \hrule height 0.5pt 
      \kern0.4ex
      \hbox{%
        \kern-0.05em
        \ensuremath{#1}%
        \kern-0.00em
      }%
    }%
  }%
}

\newcommand{\mF}{\bm{F}}

\newcommand{\mG}{\bm{G}}

\newcommand{\mU}{\bm{U}}

\newcommand{\mD}{\bm{D}}
\newcommand{\bmF}{\bm{\mathcal{F}}}
\newcommand{\bmG}{\bm{\mathcal{G}}}

\newcommand{\dt}{\Delta t}
\newcommand{\dx}{\Delta x}
\newcommand{\dy}{\Delta y}

\newcommand{\hf}{{\frac{1}{2}}}

\newcommand{\jph}{{j+\frac{1}{2}}}
\newcommand{\jmh}{{j-\frac{1}{2}}}
\newcommand{\kph}{{k+\frac{1}{2}}}
\newcommand{\kmh}{{k-\frac{1}{2}}}

\title{Adaptive High-Order A-WENO Schemes Based on a New Local Smoothness Indicator}
\author{Alina Chertock\thanks{Department of Mathematics, North Carolina State University, Raleigh, NC 27695, USA;
{\tt chertock@math.ncsu.edu}}, ~Shaoshuai Chu\thanks{Department of Mathematics, Southern University of Science and Technology, Shenzhen,
518055, China; {\tt chuss2019@mail.sustech.edu.cn}},~and Alexander Kurganov\thanks{Department of Mathematics, SUSTech International Center
for Mathematics and Guangdong Provincial Key Laboratory of Computational Science and Material Design, Southern University of Science and
Technology, Shenzhen, 518055, China; {\tt alexander@sustech.edu.cn}}}

\begin{document}

\date{}
\maketitle

\begin{abstract}
%
We develop new adaptive alternative weighted essentially non-oscillatory (A-WENO) schemes for hyperbolic systems of conservation laws. The
new schemes employ the recently proposed local characteristic decomposition based central-upwind numerical fluxes, the three-stage
third-order strong stability preserving Runge-Kutta time integrator, and the fifth-order WENO-Z interpolation. The adaptive strategy is
implemented by applying the limited interpolation only in the parts of the computational domain where the solution is identified as
``rough'' with the help of a smoothness indicator. We develop and use a new simple and robust local smoothness indicator (LSI), which is
applied to the solutions computed at each of the three stages of the ODE solver.

The new LSI and adaptive A-WENO schemes are tested on the Euler equations of gas dynamics. We implement the proposed LSI using the pressure,
which remains smooth at contact discontinuities, while our goal is to detect other ``rough'' areas and apply the limited interpolation
mostly in the neighborhoods of the shock waves. We demonstrate that the new adaptive schemes are highly accurate, non-oscillatory, and
robust. They outperform their fully limited counterparts (the A-WENO schemes with the same numerical fluxes and ODE solver but with the
WENO-Z interpolation employed everywhere) while being less computationally expensive.
\end{abstract}

\noindent
{\bf Key words:} Local smoothness indicator, scheme adaption, strong stability preserving Runge-Kutta methods, hyperbolic
systems of conservation laws, A-WENO schemes.

\noindent
{\bf AMS subject classification:} 65M06, 76M20, 76N15, 76L05, 35L65.

\section{Introduction}\label{sec1}
This paper focuses on developing high-order finite-difference methods for hyperbolic systems of conservation laws. We consider
one-dimensional (1-D),
\begin{equation}
\mU_t+\mF(\mU)_x=\bm0,
\label{1.1}
\end{equation}
and two-dimensional (2-D),
\begin{equation}
\mU_t+\mF(\mU)_x+\mG(\mU)_y=\bm0,
\label{1.1a}
\end{equation}
systems, though the proposed techniques can be directly extended to higher-dimensional cases. Here, $x$ and $y$ are spatial variables, $t$
is the time, $\mU\in\mathbb R^d$ is a vector of unknown functions, and $\mF:\mathbb R^d\to\mathbb R^d$ and $\mG:\mathbb R^d\to\mathbb R^d$
are nonlinear fluxes.

It is well-known that solutions of \eref{1.1a} may develop complicated wave structures, including shocks, rarefactions, and contact
discontinuities, even when the initial data are infinitely smooth. Therefore, it is challenging to develop highly accurate and robust
numerical methods for \eref{1.1a}. We refer the reader to various existing numerical methods, including high-order ones, e.g., the
monographs and review papers \cite{Hesthaven18,Leveque02,KLR20,Shu09,Shu20,Tor,BAF} and references therein.

Semi-discretization of \eref{1.1} and \eref{1.1a} offers one of the popular frameworks for constructing high-order finite-volume and
finite-difference schemes: the spatial derivatives are approximated using appropriate numerical fluxes. At the same time, the time evolution is
conducted with the help of a high-order and stable ODE solver. To achieve a high order of spatial accuracy, the numerical fluxes
must be evaluated using the point values of $\mU$ obtained by an appropriate piecewise polynomial reconstruction (interpolation) of the
computed solution. In order to enforce nonlinear stability, the reconstructions have to employ nonlinear limiters designed to prevent
spurious oscillations in the nonsmooth parts of the solutions. Popular finite-volume reconstructions, such as essentially non-oscillatory
(ENO) (see, e.g., \cite{Abg,Harten87,Harten87a,Shu20}) and weighted ENO (WENO) (see, e.g., \cite{Balsara20,Jiang96,Liu94,Shu09,Shu20}) ones
are highly accurate, but typically finite-volume ENO and WENO schemes are computationally expensive, especially in the multidimensional
case. More efficient implementations of ENO and WENO reconstructions can be carried out within the finite-difference framework in a
``dimension-by-dimension'' manner; see, e.g., \cite{Balsara16,Borges08,Castro11,Jiang96,Shu88,Shu89}. Unfortunately, the finite-difference
schemes, which are directly based on finite-volume reconstructions, rely on flux splittings, substantially increasing the amount of
numerical diffusion present in finite-volume ENO and WENO schemes. This drawback of finite-difference WENO schemes was overcome in
\cite{Jiang13} (also see \cite{Liu17}), where alternative WENO (A-WENO) schemes were introduced. A-WENO schemes employ standard
finite-volume numerical fluxes (without any need for flux splitting and related modifications), whose accuracy, in the context of
finite-difference schemes, is limited to the second order, while a high order is achieved using the flux Taylor expansion and high-order
WENO-Z interpolations, which were developed in \cite{Don20,Gao20,Jiang13,Liu17,Wang18}. For several recent A-WENO schemes based on
different finite-volume numerical fluxes, we refer the reader to \cite{Wang18,Wang20,WDKL}.

Even though WENO-Z interpolations are relatively computationally inexpensive and can be applied in a '`dimension-by-dimension'' manner, the
computational cost can be further reduced by avoiding the use of any nonlinear limiters in the smooth parts of the solution. In order to
achieve this goal, one needs to detect nonsmooth parts of the solution efficiently. This can be done in many ways using various existing smoothness indicators. In \cite{Berger89,Berger84}, discontinuities were detected using Richardson-type estimates of the local
truncation error of the solution. A more heuristic approach is examined in \cite{PRQ,Qui}, where the local wave strengths of the upwind
scheme were used as a measure of solution smoothness. In \cite{ABD08,Arandiga00,Arandiga99}, multiresolution coefficients of wavelets
expansions were used. In \cite{GT02, GT06}, the edges in the computed solution were detected using its Fourier coefficient. One can also
identify the ``rough'' parts of the computed solution using the numerical production of entropy (see, e.g., \cite{Pup0304,PupSem}), the
entropy residual (see, e.g., \cite{GuePas,GPP}), or the weak local residual (see, e.g., \cite{Dewar15,KK,KKP}).

In this paper, we develop a new, very simple, and robust local smoothness indicator (LSI) based on the Taylor expansion in time, applied to
the computed solutions obtained at each stage of the three-stage third-order strong stability preserving (SSP) Runge-Kutta solver; see,
e.g., \cite{Gottlieb01,Gottlieb11}. We first demonstrate that the proposed LSI can accurately detect smooth and nonsmooth solution regions.
We then apply the new LSI to design the following scheme adaption strategy in the context of the A-WENO schemes: we use the fifth-order
nonlinear WENO-Z interpolation in the detected ``rough'' parts of the computed solutions while employing a nonlimited fifth-order
interpolants in smooth areas.

The developed scheme adaption strategy is implemented using the recently proposed local characteristic decomposition based central-upwind
numerical flux from \cite{CCHKL_22} and applied to both the 1-D and 2-D Euler equations of gas dynamics, for which we design the LSI based
on the pressure rather than on the density or any other conservative variable. This choice is motivated by the results obtained in
\cite{Dewar15}, where it has been demonstrated that applying a nonlinear stabilization mechanism is crucial for the shock areas while
isolated linearly degenerate contact waves can be accurately captured using the nonlimited high-order reconstruction. We test the resulting
adaptive fifth-order A-WENO scheme on several numerical examples and demonstrate that it outperforms the corresponding fifth-order A-WENO
the scheme, which is implemented without the proposed adaptation, that is, employs the WENO-Z interpolation throughout the entire computational
domain.

The paper is organized as follows. In \S\ref{sec2}, we briefly describe the proposed 1-D and 2-D fifth-order A-WENO schemes. In
\S\ref{sec3}, we introduce the new LSI and then illustrate its performance on the Sod shock-tube problem for the 1-D Euler equations of
gas dynamics. In \S\ref{sec4f}, we describe 1-D and 2-D scheme adaption strategies based on the proposed LSI. In \S\ref{sec4}, we present a
number of the 1-D and 2-D numerical results to demonstrate the performance of the proposed adaptive A-WENO schemes and compare it with the
fully limited A-WENO schemes. Finally, in \S\ref{sec6}, we give concluding remarks.

\section{Fifth-Order A-WENO Schemes}\label{sec2}
In this section, we describe the fifth-order finite-difference A-WENO schemes introduced in \cite{Jiang13} (see also
\cite{Liu17,Wang18,Wang20,WDKL}).

\subsection{1-D A-WENO Schemes}\label{sec21}
We first consider the 1-D system \eref{1.1} and assume that the computational domain is covered with uniform cells $C_j:=[x_\jmh,x_\jph]$
of size $x_\jph-x_\jmh\equiv\dx$ centered at $x_j=(x_\jmh+x_\jph)/2$. We suppose that at a certain time $t\ge0$, the point values of the
computed solution, $\mU_j(t)$, are available, and in what follows, we will suppress the time-dependence of all of the indexed quantities for
the sake of brevity.

Following \cite{Jiang13}, $\mU_j$ are evolved in time by numerically solving the following system of ODEs:
\begin{equation}
\frac{{\rm d}\mU_j}{{\rm d}t}=-\frac{\bm{{\cal F}}_\jph-\bm{{\cal F}}_\jmh}{\dx},
\label{2.1}
\end{equation}
where $\bm{{\cal F}}_\jph$ is the fifth-order accurate numerical flux defined by
\begin{equation}
\bm{{\cal F}}_\jph\big(\bm U_\jph^-,\bm U_\jph^+\big)=\bm{{\cal F}}_\jph^{\,\rm FV}\big(\bm U_\jph^-,\bm U_\jph^+\big)-
\frac{1}{24}(\dx)^2(\mF_{xx})_\jph+\frac{7}{5760}(\dx)^4(\mF_{xxxx})_\jph.
\label{2.2f}
\end{equation}
Here, $\bm{{\cal F}}_\jph^{\,\rm FV}$ is a finite-volume numerical flux, and $({\mF_{xx}})_\jph$ and $({\mF_{xxxx}})_\jph$ are the
higher-order correction terms computed by the fourth- and second-order accurate finite differences, respectively:
\begin{equation*}
\begin{aligned}
&(\mF_{xx})_\jph=\frac{1}{48(\dx)^2}\Big[-5\mF_{j-2}+39\mF_{j-1}-34\mF_j-34\mF_{j+1}+39\mF_{j+2}-5\mF_{j+3}\Big],\\
&(\mF_{xxxx})_\jph=\frac{1}{2(\dx)^4}\Big[\mF_{j-2}-3\mF_{j-1}+2\mF_j+2\mF_{j+1}-3\mF_{j+2}+\mF_{j+3}\Big],
\end{aligned}
\end{equation*}
where $\mF_j:=\mF(\mU_j)$.

In the numerical experiments reported in \S\ref{sec41}, we have used a recently proposed local characteristics decomposition (LCD) based
central-upwind (CU) numerical flux from \cite{CCHKL_22}, which reads as
\begin{equation}
\bm{{\cal F}}_\jph^{\,\rm FV}\big(\bm U_\jph^-,\bm U_\jph^+\big)=\frac{\mF_j+\mF_{j+1}}{2}+\bm D_\jph\big(\bm U_\jph^-,\bm U_\jph^+\big),
\label{2.2}
\end{equation}
where $\bm D_\jph$ is the following numerical diffusion term:
\begin{equation}
\begin{aligned}
\bm D_\jph\big(\bm U_\jph^-,\bm U_\jph^+\big)&=R_\jph P_\jph R^{-1}_\jph\left[\mF(\bm U^-_\jph)-\frac{\mF_j+\mF_{j+1}}{2}\right]\\
&\hspace*{-1.5cm}+R_\jph M_\jph R^{-1}_\jph\left[\mF(\bm U^+_\jph)-\frac{\mF_j+\mF_{j+1}}{2}\right]
+R_\jph Q_\jph R^{-1}_\jph\left(\bm U^+_\jph-\bm U^-_\jph\right).
\end{aligned}
\label{2.3}
\end{equation}
Here, $R_\jph$ is the matrix used for the LCD in the neighborhood of $x=x_\jph$ (see Appendix \ref{appb}),
\begin{equation*}
\begin{aligned}
&P_\jph={\rm diag}\big((P_1)_\jph,\ldots,(P_d)_\jph\big),\quad M_\jph={\rm diag}\big((M_1)_\jph,\ldots,(M_d)_\jph\big),\\
&Q_\jph={\rm diag}\big((Q_1)_\jph,\ldots,(Q_d)_\jph\big)
\end{aligned}
\end{equation*}
with
$$
\begin{aligned}
&\hspace*{-0.7cm}\big((P_i)_\jph,(M_i)_\jph,(Q_i)_\jph\big)\\
&=\left\{\begin{aligned}
&\frac{1}{(\lambda^+_i)_\jph-(\lambda^-_i)_\jph}\big((\lambda^+_i)_\jph,-(\lambda^-_i)_\jph,(\lambda^+_i)_\jph(\lambda^-_i)_\jph\big)&&
\mbox{if}~(\lambda^+_i)_\jph-(\lambda^-_i)_\jph> \varepsilon,\\
&0&&\mbox{otherwise},
\end{aligned}\right.\nonumber
\end{aligned}
$$
where the one-sided local characteristic speeds,
\begin{equation}
\begin{aligned}
(\lambda^+_i)_\jph&=\max\left\{\lambda_i\big(A(\mU^-_\jph)\big),\,\lambda_i\big(A(\mU^+_\jph)\big),\, 0\right\},\\
(\lambda^-_i)_\jph&=\min\left\{\lambda_i\big(A(\mU^-_\jph)\big),\,\lambda_i\big(A(\mU^+_\jph)\big),\, 0\right\},
\end{aligned}\qquad i=1,\ldots,d,
\label{2.4}
\end{equation}
are computed using the eigenvalues $\lambda_1(A)\le\ldots\le\lambda_d(A)$ of the Jacobian $A=\frac{\partial\mF}{\partial\mU}$, and
$\varepsilon$ is a very small desingularization constant, taken $\varepsilon=10^{-10}$ in all of the numerical examples reported in
\S\ref{sec4}.

In \eref{2.2}--\eref{2.4}, $\mU^\pm_\jph$ are the right/left-sided values of $\mU$ at the cell interface $x=x_\jph$. In order to ensure the
desired fifth order of accuracy, one needs to use a fifth order accurate approximation of the point values $\mU^\pm_\jph$. It is also
important to guarantee that the resulting scheme is (essentially) non-oscillatory. This can be done by implementing a certain nonlinear
limiting procedure like the fifth-order WENO-Z interpolation from \cite{Don20,Gao20,Jiang13,Liu17,Wang18} (see Appendix \ref{appa}) applied
to the local characteristic variables (see Appendix \ref{appb}), or a certain adaption strategy like the one we will introduce in
\S\ref{sec4f}.

\subsection{2-D A-WENO Schemes}\label{sec22}
We now consider the 2-D system \eref{1.1a} and describe 2-D fifth-order A-WENO schemes.

Assume that the computational domain is covered with uniform cells $C_{j,\,k}:=[x_\jmh,x_\jph]\times[y_\kmh,y_\kph]$ centered at
$(x_j,y_k)=\big((x_\jmh+x_\jph)/2,(y_\kph+y_\kmh)/2\big)$ with $x_\jph-x_\jmh\equiv\dx$ and $y_\kph-y_\kmh\equiv\dy$ for all $j,k$. We also
assume that the computed point values $\mU_{j,k}\approx\mU(x_j,y_k,t)$ are available at a certain time level $t$. We then evolve $\mU_{j,k}$
in time by numerically solving the following system of ODEs:
\begin{equation}
\begin{aligned}
\frac{{\rm d}\mU_{j,k}}{{\rm d}t}=-\frac{\bmF_{\jph,k}-\bmF_{\jmh,k}}{\dx}-\frac{\bmG_{j,\kph}-\bmG_{j,\kmh}}{\dy},
\end{aligned}
\label{2.6f}
\end{equation}
where
$\bm{{\cal F}}_{\jph,k}$  and $\bm{{\cal G}}_{j,\kph}$ are the fifth-order accurate numerical fluxes defined by
\begin{equation}
\begin{aligned}
&\bm{{\cal F}}_{\jph,k}=\bm{{\cal F}}_{\jph,k}^{\,\rm FV}\big(\mU^{\rm E}_{j,k},\mU^{\rm W}_{j+1,k}\big)-
\frac{1}{24}(\dx)^2(\mF_{xx})_{\jph,k}+\frac{7}{5760}(\dx)^4(\mF_{xxxx})_{\jph,k},\\
&\bm{{\cal G}}_{j,\kph}=\bm{{\cal G}}_{j,\kph}^{\,\rm FV}\big(\mU^{\rm N}_{j,k},\mU^{\rm S}_{j,k+1}\big)-
\frac{1}{24}(\dy)^2(\mG_{yy})_{j,\kph}+\frac{7}{5760}(\dy)^4(\mG_{yyyy})_{j,\kph}.
\end{aligned}
\label{2.7f}
\end{equation}
Here, $\bm{{\cal F}}_{\jph,k}^{\,\rm FV}$ and $\bm{{\cal G}}_{j,\kph}^{\,\rm FV}$ are finite-volume fluxes, and $(\mF_{xx})_{\jph,k}$,
$(\mG_{yy})_{j,\kph}$, $(\mF_{xxxx})_{\jph,k}$, and $(\mG_{yyyy})_{j,\kph}$ are the higher-order correction terms computed by the fourth-
and second-order accurate finite differences, respectively:
\begin{equation*}
\begin{aligned}
&(\mF_{xx})_{\jph,k}=\frac{1}{48(\dx)^2}\left(-5\mF_{j-2,k}+39\mF_{j-1,k}-34\mF_{j,k}-34\mF_{j+1,k}+39\mF_{j+2,k}-5\mF_{j+3,k}\right),\\
&(\mF_{xxxx})_{\jph,k}=\frac{1}{2(\dx)^4}\left(\mF_{j-2,k}-3\mF_{j-1,k}+2\mF_{j,k}+2\mF_{j+1,k}-3\mF_{j+2,k}+\mF_{j+3,k}\right),\\
&(\mG_{yy})_{j,\kph}=\frac{1}{48(\dy)^2}\left(-5\mG_{j,k-2}+39\mG_{j,k-1}-34\mG_{j,k}-34\mG_{j,k+1}+39\mG_{j,k+2}-5\mG_{j,k+3}\right),\\
&(\mG_{yyyy})_{j,\kph}=\frac{1}{2(\dy)^4}\left(\mG_{j,k-2}-3\mG_{j,k-1}+2\mG_{j,k}+2\mG_{j,k+1}-3\mG_{j,k+2}+\mG_{j,k+3}\right),
\end{aligned}
\end{equation*}
where $\mF_{j,k}:=\mF(\mU_{j,k})$ and $\mG_{j,k}:=\mG(\mU_{j,k})$.

In the numerical experiments reported in \S\ref{sec42}, we have used the 2-D LCD-based CU numerical fluxes from \cite{CCHKL_22}:
\begin{equation}
\begin{aligned}
&\bmF^{\rm FV}_{\jph,k}\big(\mU^{\rm E}_{j,k},\mU^{\rm W}_{j+1,k}\big)=\frac{\mF_{j,k}+\mF_{j+1,k}}{2}+
\mD_{\jph,k}\big(\mU^{\rm E}_{j,k},\mU^{\rm W}_{j+1,k}\big),\\
&\bmG^{\rm FV}_{j,\kph}\big(\mU^{\rm N}_{j,k},\mU^{\rm S}_{j,k+1}\big)=
\frac{\mG_{j,k}+\mG_{j,k+1}}{2}+\mD_{j,\kph}\big(\mU^{\rm N}_{j,k},\mU^{\rm S}_{j,k+1}\big),
\end{aligned}
\label{2.5}
\end{equation}
where $\mD_{\jph,k}$ and $\mD_{j,\kph}$ are numerical diffusion terms defined by
\begin{equation}
\hspace*{-0.2cm}\begin{aligned}
&\mD_{\jph,k}\big(\mU^{\rm E}_{j,k},\mU^{\rm W}_{j+1,k}\big)=
R_{\jph,k}P_{\jph,k}R^{-1}_{\jph,k}\left[\mF\big(\bm U^{\rm E}_{j,k}\big)-\frac{\mF_{j,k}+\mF_{j+1,k}}{2}\right]\\
&+R_{\jph,k}M_{\jph,k}R^{-1}_{\jph,k}\left[\mF\big(\bm U^{\rm W}_{j+1,k}\big)-
\frac{\mF_{j,k}+\mF_{j+1,k}}{2}\right]+R_{\jph,k}Q_{\jph,k}R^{-1}_{\jph,k}\left(\mU^{\rm W}_{j+1,k}-\mU^{\rm E}_{j,k}\right),\\
&\mD_{j,\kph}\big(\mU^{\rm N}_{j,k},\mU^{\rm S}_{j,k+1}\big)=
R_{j,\kph}P_{j,\kph}R^{-1}_{j,\kph}\left[\mG\big(\bm U^{\rm N}_{j,k}\big)-\frac{\mG_{j,k}+\mG_{j,k+1}}{2}\right]\\
&+R_{j,\kph}M_{j,\kph}R^{-1}_{j,\kph}\left[\mG\big(\bm U^{\rm S}_{j,k+1}\big)-\frac{\mG_{j,k}+\mG_{j,k+1}}{2}\right]+
R_{j,\kph}Q_{j,\kph}R^{-1}_{j,\kph}\left(\mU^{\rm S}_{j,k+1}-\mU^{\rm N}_{j,k}\right).
\end{aligned}
\label{2.6}
\end{equation}
The matrices $R_{\jph,k}$ and $R_{j,\kph}$ are used for the LCD in the neighborhoods of $(x,y)=(x_\jph,y_k)$ and $(x,y)=(x_j, y_\kph)$,
respectively, and
\begin{equation*}
\begin{aligned}
&P_{\jph,k}={\rm diag}\big(\big(P_1\big)_{\jph,k},\ldots,\big(P_d\big)_{\jph,k}\big),&&
P_{j,\kph}={\rm diag}\big(\big(P_1\big)_{j,\kph},\ldots,\big(P_d\big)_{j,\kph}\big),\\
&M_{\jph,k}={\rm diag}\big(\big(M_1\big)_{\jph,k},\ldots,\big(M_d\big)_{\jph,k}\big),&&
M_{j,\kph}={\rm diag}\big(\big(M_1\big)_{j,\kph},\ldots,\big(M_d\big)_{j,\kph}\big),\\
&Q_{\jph,k}={\rm diag}\big(\big(Q_1\big)_{\jph,k},\ldots,\big(Q_d\big)_{\jph,k}\big),&&
Q_{j,\kph}={\rm diag}\big(\big(Q_1\big)_{j,\kph},\ldots,\big(Q_d\big)_{j,\kph}\big),
\end{aligned}
\end{equation*}
with
$$
\begin{aligned}
&\hspace*{-1.0cm}\big((P_i)_{\jph,k},(M_i)_{\jph,k},(Q_i)_{\jph,k}\big)\\
&=\left\{\begin{aligned}
&\frac{1}{\Delta(\lambda_i)_{\jph,k}}
\left((\lambda^+_i)_{\jph,k},-(\lambda^-_i)_{\jph,k},(\lambda^+_i)_{\jph,k}(\lambda^-_i)_{\jph,k}\right)
&&\mbox{if}~\Delta(\lambda_i)_{\jph,k}>\varepsilon,\\
&0&&\mbox{otherwise},
\end{aligned}\right.\\
&\hspace*{-1.0cm}\big((P_i)_{j,\kph},(M_i)_{j,\kph},(Q_i)_{j,\kph}\big)\\
&=\left\{\begin{aligned}
&\frac{1}{\Delta(\mu_i)_{j,\kph}}\left((\mu^+_i)_{j,\kph},-(\mu^-_i)_{j,\kph},(\mu^+_i)_{j,\kph}(\mu^-_i)_{j,\kph}\right)
&&\mbox{if}~\Delta(\mu_i)_{j,\kph}>\varepsilon,\\
&0&&\mbox{otherwise}.
\end{aligned}\right.
\end{aligned}
$$
Here, $\Delta(\lambda_i)_{\jph,k}:=(\lambda^+_i)_{\jph,k}-(\lambda^-_i)_{\jph,k}$,
$\,\Delta(\mu_i)_{j,\kph}:=(\mu^+_i)_{j,\kph}-(\mu^-_i)_{j,\kph}$, and
\begin{equation}
\begin{aligned}
(\lambda^+_i)_{\jph,k}=\max\left\{\lambda_i\big(A(\mU^{\rm E}_{j,k})\big),\,\lambda_i\big(A(\mU^{\rm W}_{j+1,k})\big),\,0\right\},\\
(\lambda^-_i)_{\jph,k}=\min\left\{\lambda_i\big(A(\mU^{\rm E}_{j,k})\big),\,\lambda_i\big(A(\mU^{\rm W}_{j+1,k})\big),\,0\right\},\\
(\mu^+_i)_{j,\kph}=\max\left\{\mu_i\big(B(\mU^{\rm N}_{j,k})\big),\,\mu_i\big(B(\mU^{\rm S}_{j,k+1})\big),\,0\right\},\\
(\mu^-_i)_{j,\kph}=\min\left\{\mu_i\big(B(\mU^{\rm N}_{j,k})\big),\,\mu_i\big(B(\mU^{\rm S}_{j,k+1})\big),\,0\right\},
\end{aligned}
\label{2.7}
\end{equation}
where $\lambda_1(A)\le\ldots\le\lambda_d(A)$ and $\mu_1(B)\le\ldots\le\mu_d(B)$ are the eigenvalues of the Jacobians
$A=\frac{\partial\mF}{\partial\mU}$ and $B=\frac{\partial\mG}{\partial\mU}$, respectively.

In \eref{2.5}--\eref{2.7}, $\mU^{\rm E}_{j,k}$, $\mU^{\rm W}_{j+1,k}$ and $\mU^{\rm N}_{j,k}$, $\mU^{\rm S}_{j,k+1}$ are the one-sided
values of $\mU$ at the cell interfaces $(x,y)=(x_\jph\pm0,y_k)$ and $(x,y)=(x_j,y_\kph\pm0)$, respectively. In order to achieve fifth-order
accuracy, $\mU^{\rm E(W)}_{j,k}$ and $\mU^{\rm N(S)}_{j,k}$ are, as in the 1-D case, approximated either using the fifth-order WENO-Z
interpolant applied to the local characteristic variables in the $x$- and $y$-directions, respectively, or with the help of the adaptive
strategy, which we will introduce in \S\ref{sec4f}.

\section{A New Local Smoothness Indicator (LSI)}\label{sec3}
In this section, we introduce a very simple LSI, which we will later use as a base for a scheme adaption strategy.

We first consider a function $\psi(\cdot,t)$ and introduce the following quantity:
\begin{equation}
D^\psi(\cdot,t-\tau):=\Big|\frac{\psi(\cdot,t)+\psi(\cdot,t-2\tau)}{2}-\psi(\cdot,t-\tau)\Big|,
\label{3.1f}
\end{equation}
where $\tau>0$ and $\cdot$ stand for a certain spatial coordinate. If $\psi$ is smooth, then one can use the Taylor expansion about the
point $(\cdot,t-\tau)$ to obtain
\begin{equation}
D^\psi(\cdot,t-\tau)=\frac{\tau^2}{2}\psi_{tt}(\cdot,t-\tau)+{\cal O}(\tau^4).
\label{3.2f}
\end{equation}
This suggests that for piecewise smooth $\psi$ the magnitude of $D^\psi$ is proportional to $\tau^2$ in the areas where $\psi$ is smooth and
is ${\cal O}(1)$ elsewhere.

In order to design an LSI based on \eref{3.1f}, we proceed as follows. We begin with the 1-D case, denote by $\mU(t):=\{\mU_j(t)\}$, and let
${\cal L}[\mU(t)]$ be the nonlinear operator representing the right-hand side (RHS) of \eref{2.1}. Assuming that the computed solution is
available at a certain discrete time level $t=t^n$, we evolve it to the next time level $t^{n+1}:=t^n+\dt^n$ by numerically integrating the
ODE system \eref{2.1} using the three-stage third-order SSP Runge-Kutta method, which reads as (see \cite{Gottlieb01,Gottlieb11})
\begin{equation}
\begin{aligned}
&\mU^{\rm I}(t^{n+1})=\mU(t^n)+\dt^n{\cal L}[\mU(t^n)],\\
&\mU^{\rm II}\big(t^{n+\hf}\big)=\frac{3}{4}\,\mU(t^n)+\frac{1}{4}\left(\mU^{\rm I}(t^{n+1})+\dt^n{\cal L}[\mU^{\rm I}(t^{n+1})]\right),\\
&\mU(t^{n+1})=\frac{1}{3}\,\mU(t^n)+\frac{2}{3}\left(\mU^{\rm II}\big(t^{n+\hf}\big)+\dt^n{\cal L}\big[\mU^{\rm II}\big(t^{n+\hf}\big)\big]
\right),
\end{aligned}
\label{3.1}
\end{equation}
where $t^{n+\hf}:=t^n+\dt^n/2$ and $\mU^{\rm I}$ and $\mU^{\rm II}$ are the intermediate stage solutions, which are lower-order
approximations of $\mU$ at time levels $t^{n+1}$ and $t^{n+\hf}$, respectively. In \eref{3.1}, the time step $\dt^n$ is selected based on
the following CFL-based stability restriction:
\begin{equation}
\dt^n\le\frac{\dx}{2a},\quad a:=\max_j\left\{\max\big((\lambda_d^+)_\jph,-(\lambda_1^-)_\jph\big)\right\}.
\label{3.4ff}
\end{equation}

Next, we introduce quantities $\psi_j(t^n):=\psi(\mU_j(t^n))$ and
$\psi_j^{\rm II}\big(t^{n-\hf}\big):=\psi\big(\mU_j^{\rm II}\big(t^{n-\hf}\big)\big)$, and the corresponding LSI based on \eref{3.1f} with
$\tau=\dt^{n-1}/2$:
\begin{equation}
D_j^\psi\big(t^{n-\hf}\big)=\Big|\frac{\psi_j(t^n)+\psi_j(t^{n-1})}{2}-\psi_j^{\rm II}\big(t^{n-\hf}\big)\Big|.
\label{3.4f}
\end{equation}

Similarly, the 2-D LSI is given by
\begin{equation}
D_{j,k}^\psi\big(t^{n-\hf}\big)=\Big|\frac{\psi_{j,k}(t^n)+\psi_{j,k}(t^{n-1})}{2}-\psi_{j,k}^{\rm II}\big(t^{n-\hf}\big)\Big|,
\label{3.5f}
\end{equation}
where $\psi_{j,k}(t^n):=\psi(\mU_{j,k}(t^n))$ and
$\psi_{j,k}^{\rm II}\big(t^{n-\hf}\big):=\psi\big(\mU_{j,k}^{\rm II}\big(t^{n-\hf}\big)\big)$.

While the computation in \eref{3.2f} is based on the smoothness of $\psi$, the LSIs \eref{3.4f} and \eref{3.5f} can, in principle, be used
for detecting ``rough'' areas of nonsmooth computed solutions. However, before these LSIs can be used for the development of the robust
adaptation strategies, one may need to smear the introduced quantities in space by introducing
\begin{equation}
\xbar D_j^{\,\psi}\big(t^{n-\hf}\big):=\frac{1}{6}
\left[D_{j-1}^\psi\big(t^{n-\hf})+D_j^\psi\big(t^{n-\hf})+D_{j+1}^\psi\big(t^{n-\hf})\right]
\label{3.7f}
\end{equation}
and
\begin{equation*}
\begin{aligned}
\xbar D_{j,k}^{\,\psi}\big(t^{n-\hf}\big)&=\frac{1}{36}\left[D_{j-1,k-1}^\psi\big(t^{n-\hf}\big)+D_{j-1,k+1}^\psi\big(t^{n-\hf}\big)+
D_{j+1,k-1}^\psi\big(t^{n-\hf}\big)+D_{j+1,k+1}^\psi\big(t^{n-\hf}\big)\right.\\
&\left.+\,4\left(D_{j-1,k}^\psi\big(t^{n-\hf}\big)+D_{j,k-1}^\psi\big(t^{n-\hf}\big)+D_{j,k+1}^\psi\big(t^{n-\hf}\big)+
D_{j+1,k}^\psi\big(t^{n-\hf}\big)\right)+16D_{j,k}^\psi\big(t^{n-\hf}\big)\right]
\end{aligned}
\end{equation*}
in the 1-D and 2-D cases, respectively.

In order to verify the plausibility of possible adaptation strategies based on the introduced LSI, we measure its size in the following
a numerical example in which we compute the solution of a benchmark using the fully limited A-WENO scheme that employs the WENO-Z interpolant
throughout the entire computational domain (see \S\ref{sec21}).

\paragraph{Example---Sod Shock-Tube Problem for Euler Equations of Gas Dynamics.}
We consider the 1-D Euler equations of gas dynamics, which reads as
\begin{equation}
\begin{aligned}
&\rho_t+(\rho u)_x=0,\\
&(\rho u)_t+(\rho u^2+p)_x=0,\\
&E_t+\left[u(E+p)\right]_x=0,
\end{aligned}
\label{3.6}
\end{equation}
where $\rho$, $u$, $p$, and $E$ are the density, velocity, pressure, and total energy, respectively. The system \eref{3.6} is completed
through the following equations of state:
\begin{equation}
p=(\gamma-1)\Big[E-\hf\rho u^2\Big],
\label{3.7}
\end{equation}
where the parameter $\gamma$ represents the specific heat ratio (we take $\gamma=1.4$). We consider the following initial conditions
\cite{Sod78}:
\begin{equation}
(\rho,u,p)(x,0)=\begin{cases}(1,0,1.0),&x<0.5,\\(0.125,0,0.1),&x>0.5,\end{cases}
\label{3.7a}
\end{equation}
prescribed in the interval $[0,1]$ subject to the free boundary conditions.

\begin{figure}[ht!]
\centerline{\includegraphics[trim=0.8cm 0.3cm 1.3cm 0.cm, clip, width=6.0cm]{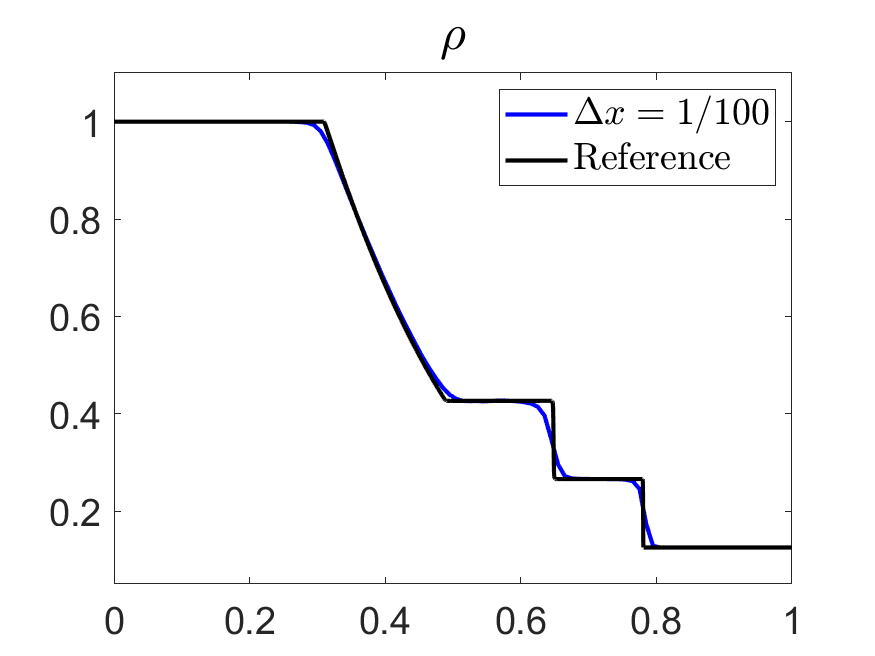}\hspace{1cm}
            \includegraphics[trim=0.8cm 0.3cm 1.3cm 0.cm, clip, width=6.0cm]{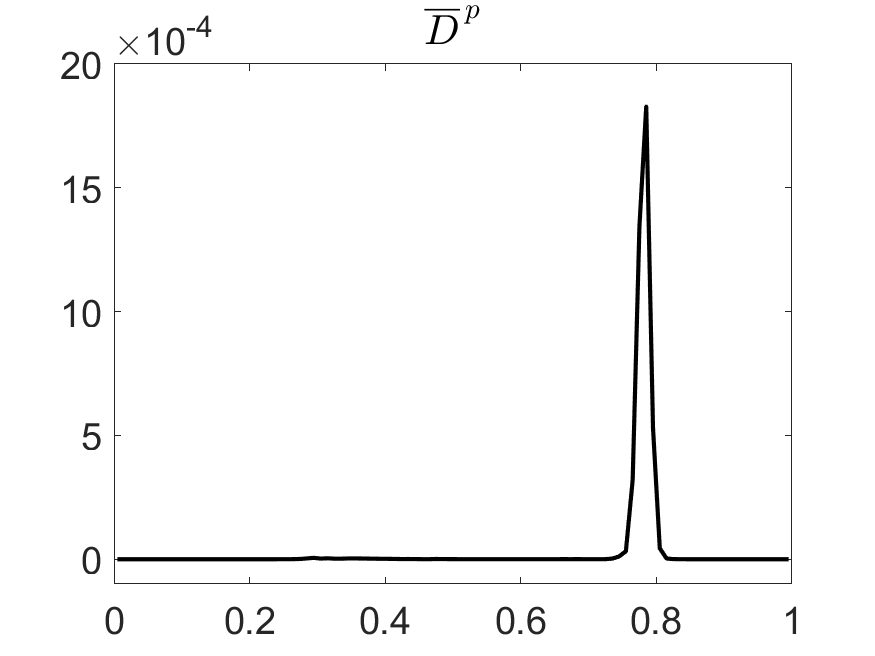}}
\caption{\sf Density (left) and the corresponding values of the pressure-based LSI (right).\label{fig1}}
\end{figure}

We compute the numerical solution by the A-WENO scheme introduced in \S\ref{sec21} until the final time $t=0.16$ on a uniform mesh with
$\dx=1/100$. In Figure \ref{fig1}, we plot the obtained density together with the reference solution computed on a much finer uniform mesh
with $\dx=1/4000$ and the pressure-based LSI ($\xbar D_j^{\,p}$ defined in \eref{3.7f} with $\psi=p$) computed at the final time step. As
one can see, the LSI can detect the shock wave's location and indicate the area of a rarefaction corner. At the same time,
the LSI values in the contact discontinuity neighborhood are very small. If, however, one is interested in identifying contact
discontinuities as well, one can use the density-based LSI ($\xbar D_j^{\,\rho}$ defined in \eref{3.7f} with $\psi=\rho$).

In order to investigate the plausibility of the LSI-based adaptive strategies, we compute the numerical solutions on a sequence
of uniform meshes with $\dx=1/200$, 1/400, 1/800, 1/1600, 1/3200, and 1/6400, and measure the asymptotic behavior of the LSI in different
parts of the computational domain. The obtained results are reported in Table \ref{table1}, where one can observe quite significant
differences in the order of magnitude of the LSI. For example, on the mesh with $\dx=1/200$, the local maxima of the $\xbar D^{\,p}$ are
$\sim10^{-6}$ near the rarefaction corner and in a smooth region within the rarefaction wave, $\sim10^{-7}$ near the contact discontinuity,
and $\sim10^{-3}$ at the shock (the last local maximum is, in fact, the global maximum of $\xbar D^{\,p}$). One can also see that away
from the shock, the LSI decays when the mesh is refined. The rate of decay is second-order in the smooth region and about the first-order
near the rarefaction corner. In the contact wave area, the LSI is very small, and when the coarse mesh is refined, the LSI decays very rapidly
there. At the same time, near the shock, the size of $\xbar D^{\,p}$ is practically independent of the mesh size, as expected. This suggests
that the proposed LSI can be used as an efficient and accurate tool to detect shocks and other ``rough'' parts of the computed solution
except for the isolated contact waves, which can be treated in the same way as smooth parts of the computed solution; see the description of
the adaption strategy we propose in the next section.
\begin{table}[ht!]
\centering
\begin{tabular}{c cc cc}
\hline
$\dx$&$\max\limits_{a\le x_j\le b}\xbar D_j^{\,p}$&Rate&$\max\limits_{a\le x_j\le b}\xbar D_j^{\,p}$&Rate\\
\hline
\multicolumn{2}{l}{Rarefaction corner, $a=0.25, b=0.35$} & &\multicolumn{2}{l}{Smooth subregion, $a=0.35, b=0.45$}\\
\hline
$1/100$ &7.94e-06&--  &3.87e-06&--  \\
$1/200$ &2.86e-06&1.47&1.09e-06&1.83 \\
$1/400$ &1.28e-06&1.17&2.72e-07&2.01\\
$1/800$ &4.81e-07&1.41&6.75e-08&2.01\\
$1/1600$&2.27e-07&1.08&1.66e-08&2.02\\
$1/3200$&9.76e-08&1.22&4.14e-09&2.01\\
$1/6400$&4.43e-08&1.14&1.03e-09&2.00\\
\hline
\multicolumn{2}{l}{Contact wave, $a=0.6, b=0.7$} & &\multicolumn{2}{l}{Everywhere (Shock), $a=0, b=1$}\\
\hline
$1/100$ &6.59e-07&--  &2.58e-03&--\\
$1/200$ &1.33e-07&2.31&2.36e-03&0.13\\
$1/400$ &1.54e-09&6.43&1.82e-03&0.38\\
$1/800$ &6.55e-10&1.23&2.37e-03&-0.38\\
$1/1600$&3.50e-10&0.90&5.76e-04&2.04\\
$1/3200$&2.00e-10&0.80&2.00e-03&-1.79\\
$1/6400$&8.28e-11&1.27&2.39e-03&-0.26\\
\hline
\end{tabular}
\caption{\sf Local and global maxima of $\xbar D^{\,p}$ and the corresponding rates of change.\label{table1}}
\end{table}

\section{Scheme Adaption}\label{sec4f}
In this section, we develop a scheme adaption strategy based on the LSIs from \S\ref{sec3} and the A-WENO schemes described in \S\ref{sec2}.
This will lead to new adaptive A-WENO schemes, in which the WENO-Z interpolation will only be used in the ``rough'' areas indicated by the
LSI.

\paragraph{One-Dimensional Algorithm.} Assume that $\mU_j(t^n)$, $\mU_j(t^{n-1})$, and $\mU_j^{\rm II}\big(t^{n-\hf}\big)$ are available for
all $j$. We then compute the LSI values given by \eref{3.7f} and identify the ``rough'' areas as follows. We first find all of the points
$x=x_j$ at which
\begin{equation}
\xbar D^{\,\psi}_j\big(t^{n-\hf}\big)>\texttt{C}(\dt^{n-1})^\frac{3}{2},
\label{4.1}
\end{equation}
where $\texttt{C}$ is a positive tunable constant to be selected for each problem at hand, and presume that the solution at time
$t=t^{n-\hf}$ is ``rough'' there. Due to the finite speed of propagation and the CFL condition \eref{3.4ff}, one may presume that the
solution at the time interval $[t^n,t^{n+1}]$ (that is, at all of the three Runge-Kutta stages \eref{3.1}) is ``rough'' at the nearby points
$x_{j\pm\hf}$ and $x_{j\pm\frac{3}{2}}$.

After identifying each of the points $x_\jph$ as either ``rough'' or ``smooth'', we compute either nonlimited $\breve\mU_\jph^\pm$ or
limited $\widetilde\mU_\jph^\pm$ point values there, and then evaluate the finite-volume numerical fluxes needed in \eref{2.1}--\eref{2.2f}
(and hence in \eref{3.1}) by
$$
\bm{{\cal F}}_\jph^{\,\rm FV}=\left\{\begin{aligned}
&\bm{{\cal F}}_\jph^{\,\rm FV}\big(\breve\mU_\jph^-,\breve\mU_\jph^+\big)&&\mbox{if}~x_\jph~\mbox{is ``rough''},\\
&\bm{{\cal F}}_\jph^{\,\rm FV}\big(\widetilde\mU_\jph^-,\widetilde\mU_\jph^+\big)&&\mbox{if}~x_\jph~\mbox{is ``smooth''}.
\end{aligned}\right.
$$

\paragraph{Two-Dimensional Algorithm.} An extension of the 1-D scheme adaption algorithm to the 2-D case is relatively straightforward.

The main component of the 2-D algorithm is identifying the ``rough'' parts of the solution, in which the one-sided interpolated
values are to be computed using the WENO-Z interpolant. As in the 1-D case, this is done using the LSI. Namely, we presume that the solution
at time $t=t^{n-\hf}$ is ``rough'' in all of the cells $C_{j,k}$, in which
\begin{equation}
\xbar D^{\,\psi}_{j,k}\big(t^{n-\hf}\big)>\texttt{C}(\dt^{n-1})^\frac{3}{2}.
\label{4.2}
\end{equation}
Then, due to the finite speed of propagation and the appropriate CFL condition with the CFL number 1/2, one may presume that the solution at
the time interval $[t^n,t^{n+1}]$ is ``rough'' at the nearby points $(x_{j\pm\hf},y_{k\pm1})$, $(x_{j\pm\hf},y_k)$,
$(x_{j\pm\frac{3}{2}},y_k)$ and $(x_{j\pm1},y_{k\pm\hf})$, $(x_j,y_{k\pm\hf})$, $(x_j,y_{k\pm\frac{3}{2}})$.

Equipped with the information about the ``rough'' and ``smooth'' parts of the computed solutions, we proceed with the proposed adaption
strategy and compute either nonlimited $\breve\mU_{j,k}^{\rm E(W,N,S)}$ or limited $\widetilde\mU_{j,k}^{\rm E(W,N,S)}$ point values there,
and then evaluate the finite-volume numerical fluxes needed in \eref{2.6f}--\eref{2.7f} by
$$
\begin{aligned}
\bm{{\cal F}}_{\jph,k}^{\,\rm FV}&=\left\{\begin{aligned}
&\bm{{\cal F}}_{\jph,k}^{\,\rm FV}\big(\breve\mU_{j,k}^{\,\rm E},\breve\mU_{j+1,k}^{\,\rm W}\big)&&
\mbox{if}~(x_\jph,y_k)~\mbox{is ``rough''},\\
&\bm{{\cal F}}_{\jph,k}^{\,\rm FV}\big(\widetilde\mU_{j,k}^{\,\rm E},\widetilde\mU_{j+1,k}^{\,\rm W}\big)&&
\mbox{if}~(x_\jph,y_k)~\mbox{is ``smooth''},
\end{aligned}\right.\\
\bm{{\cal G}}_{j,\kph}^{\,\rm FV}&=\left\{\begin{aligned}
&\bm{{\cal G}}_{j,\kph}^{\,\rm FV}\big(\breve\mU_{j,k}^{\,\rm N},\breve\mU_{j,k+1}^{\,\rm S}\big)&&
\mbox{if}~(x_j,y_\kph)~\mbox{is ``rough''},\\
&\bm{{\cal G}}_{j,\kph}^{\,\rm FV}\big(\widetilde\mU_{j,k}^{\,\rm N},\widetilde\mU_{j,k+1}^{\,\rm S}\big)&&
\mbox{if}~(x_j,y_\kph)~\mbox{is ``smooth''}.
\end{aligned}\right.
\end{aligned}
$$
\begin{rmk}
The fact that the constants $\texttt{C}$ in \eref{4.1} and \eref{4.2} must be tuned is a weak point of our adaption strategy. One may,
however, tune $\texttt{C}$ on a coarse mesh and then use the same value of $\texttt{C}$ on finer meshes to minimize an extra computational
cost as it was done, e.g., in \cite{Kurganov12a} in the context of an adaptive artificial viscosity method. The plausibility of this strategy
in the current scheme adaption algorithm is supported by a numerical experiment; see Example 1 in \S\ref{sec41} below.
\end{rmk}
\begin{rmk}
As no past time solution is available at the first time step, at $t^0=0$, we complete the first evolution step using a fully
limited A-WENO scheme that employs the WENO-Z interpolation throughout the entire computational domain.
\end{rmk}

\section{Numerical Examples}\label{sec4}
In this section, we test the developed adaption strategy on several numerical examples. To this end, we apply the adaptive A-WENO
schemes to several initial-boundary value problems for the 1-D and 2-D Euler equations of gas dynamics and compare their performance with
the fully limited A-WENO schemes. In the rest of this section, we will refer to the proposed adaptive A-WENO schemes as
to {\em adaptive schemes} and the fully limited A-WENO schemes as to {\em limited schemes}.

In all of the examples below, the specific heat ratio is $\gamma=1.4$ (except for Example 9, where $\gamma=5/3$), and the CFL number is
0.45.

\subsection{One-Dimensional Examples}\label{sec41}
\begin{figure}[ht!]
\centerline{\includegraphics[trim=0.5cm 0.3cm 1.3cm 0.4cm, clip, width=6.0cm]{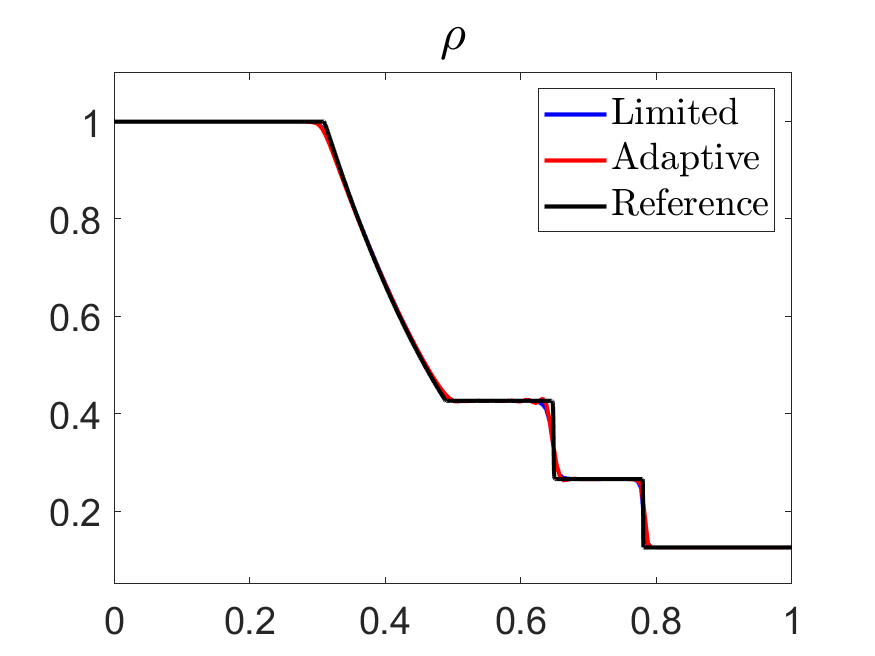}\hspace{1cm}
            \includegraphics[trim=0.5cm 0.3cm 1.3cm 0.4cm, clip, width=6.0cm]{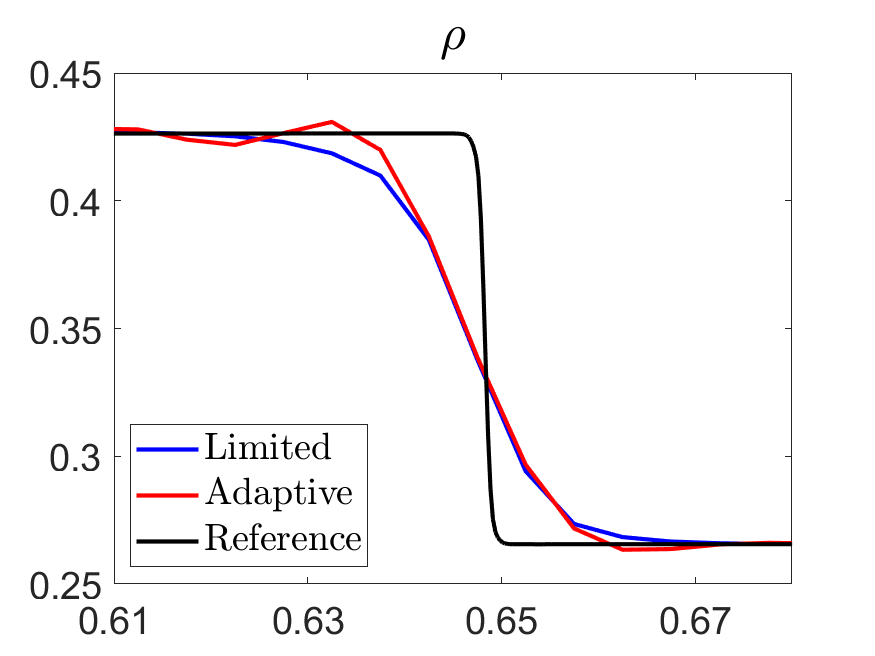}}
\caption{\sf Example 1: Density $\rho$ computed by the limited and adaptive schemes (left) and zoom at $x\in[0.61,0.68]$ (right).
\label{fig3}}
\end{figure}

\begin{figure}[ht!]
\centerline{\includegraphics[trim=0.8cm 0.3cm 1.3cm 0.1cm, clip, width=6.0cm]{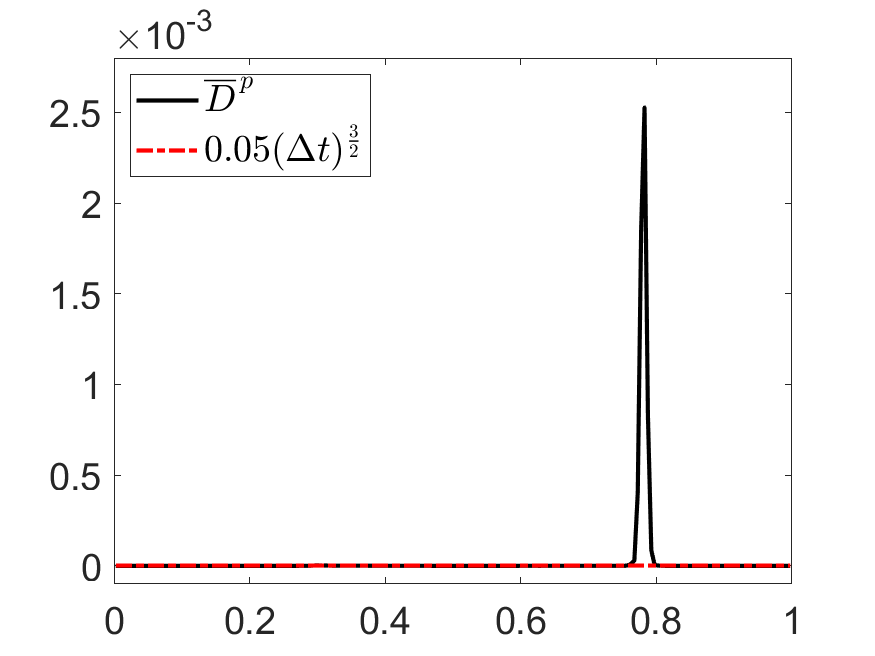}\hspace{1cm}
            \includegraphics[trim=0.8cm 0.3cm 1.3cm 0.1cm, clip, width=6.0cm]{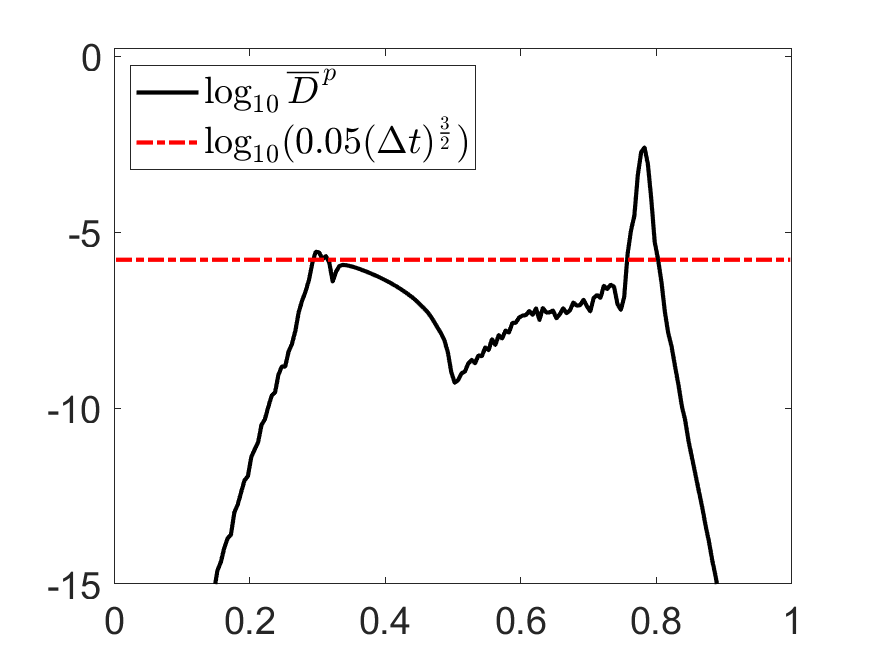}}
\caption{\sf Example 1: $\xbar D^{\,p}$ and $0.05(\dt)^\frac{3}{2}$ (left) and $\log_{10}\xbar D^{\,p}$ and
$\log_{10}(0.05(\dt)^\frac{3}{2})$ (right) for $\dx=1/200$.\label{fig4a}}
\end{figure}

\paragraph{Example 1---Sod Shock-Tube Problem.} In the first example, we once again consider \eqref{3.6}--\eqref{3.7a} subject to the free
boundary conditions and compute the numerical solution until the final time $t=0.16$ using both the limited and adaptive schemes. In this
example, we take $\texttt{C}=0.05$ while implementing the scheme adaption strategy. The obtained solutions, computed on a uniform mesh with
$\dx=1/200$ and the corresponding reference solution computed by the limited scheme on a much finer mesh with $\dx=1/4000$ are presented in
Figure \ref{fig3}. We also plot the LSI $\xbar D^{\,p}$ and $0.05(\dt)^\frac{3}{2}$ along with $\log_{10}\xbar D^{\,p}$ and
$\log_{10}(0.05(\dt)^\frac{3}{2})$, computed during the adaptive scheme evolution at the final time moment in Figure \ref{fig4a}. One can
observe that the computed LSI can capture the shock wave position accurately, and the results obtained by the adaptive scheme are sharper than those obtained by the limited scheme, even though there are small oscillations near the contact wave captured by the
adaptive scheme.

\begin{figure}[ht!]
\centerline{\includegraphics[trim=0.8cm 0.3cm 1.3cm 0.1cm, clip, width=6.0cm]{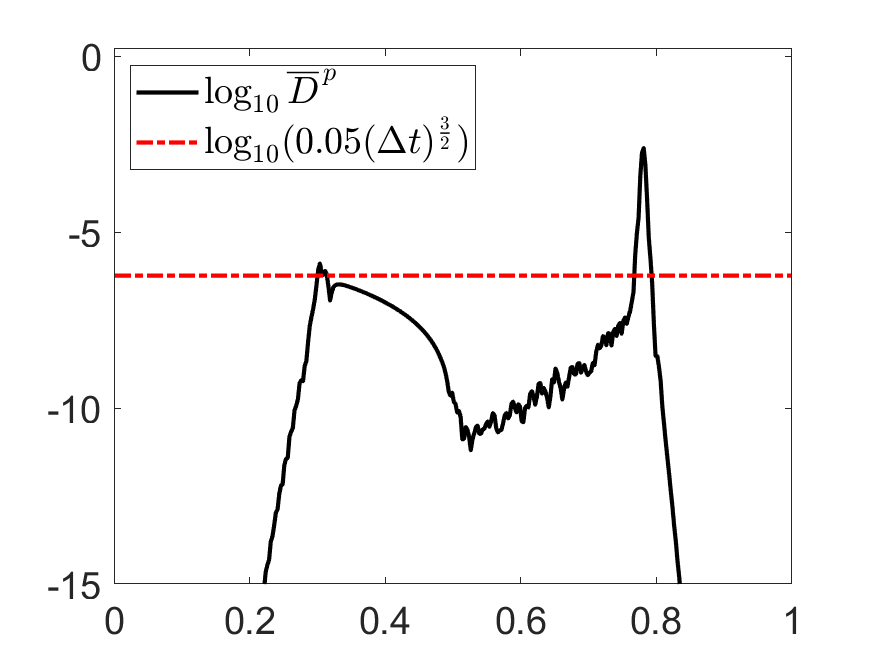}\hspace{1cm}
            \includegraphics[trim=0.8cm 0.3cm 1.3cm 0.1cm, clip, width=6.0cm]{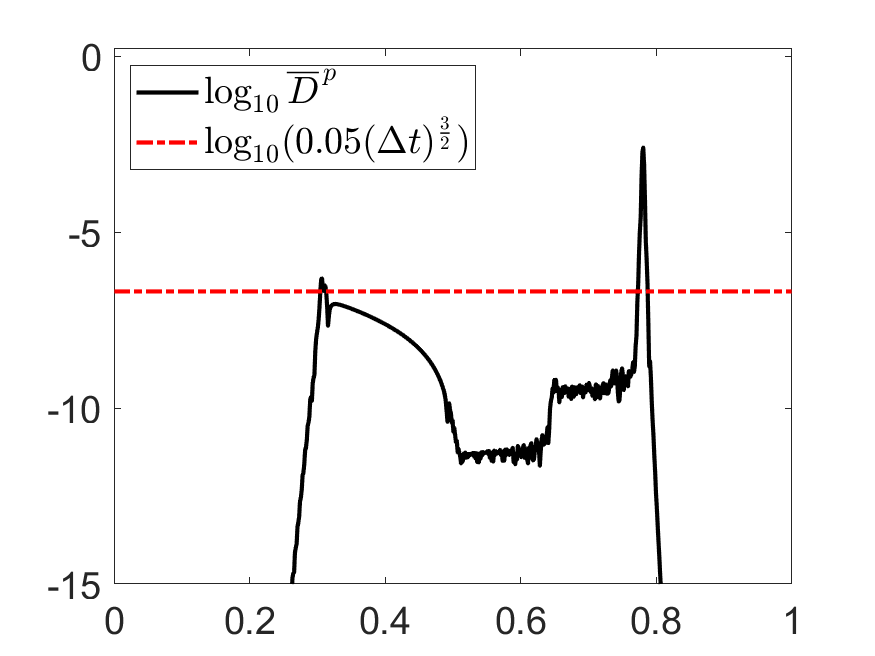}}
\caption{\sf Example 1: $\log_{10}\xbar D^{\,p}$ and $\log_{10}(0.05(\dt)^\frac{3}{2})$ for $\dx=1/400$ (left) and $1/800$ (right).
\label{fig4b}}
\end{figure}

It is also instructive to point out that choosing the adaption constant $\texttt{C}$ on a coarse mesh is a robust strategy. This is evident
from the graphs of $\log_{10}\xbar D^{\,p}$ and $\log_{10}(0.05(\dt)^\frac{3}{2})$ depicted in Figures \ref{fig4a} (right) and \ref{fig4b}.
The presented results illustrate that the ``rough'' parts of the computed solutions can be accurately identified using \eref{4.1}, while the
same constant $\texttt{C}=0.05$ is used on three different meshes.

\paragraph{Example 2---``Shock-Bubble'' Iteration Problem.} In the second example taken from \cite{KX_22}, we consider the ``shock-bubble''
interaction problem. The initial data for the 1-D Euler equations \eref{3.6}--\eref{3.7},
\begin{equation*}
(\rho, u,p)(x,0)=\begin{cases}
(13.1538,0,1),&|x|<0.25,\\
(1.3333,-0.3535,1.5),&x>0.75,\\
(1,0,1),&\mbox{otherwise,}
\end{cases}
\end{equation*}
correspond to a left-moving shock, initially located at $x=0.75$, and a ``bubble'' with a radius of 0.25, initially located at the origin.

\begin{figure}[ht!]
\centerline{\includegraphics[trim=0.8cm 0.3cm 1.2cm 0.4cm, clip, width=6.0cm]{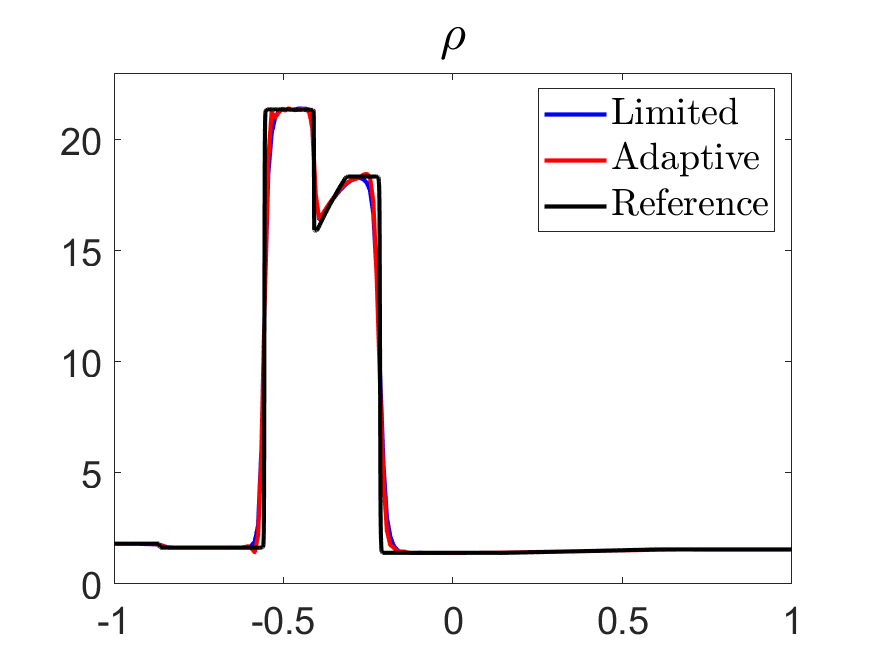}\hspace{1cm}
            \includegraphics[trim=0.8cm 0.3cm 1.2cm 0.4cm, clip, width=6.0cm]{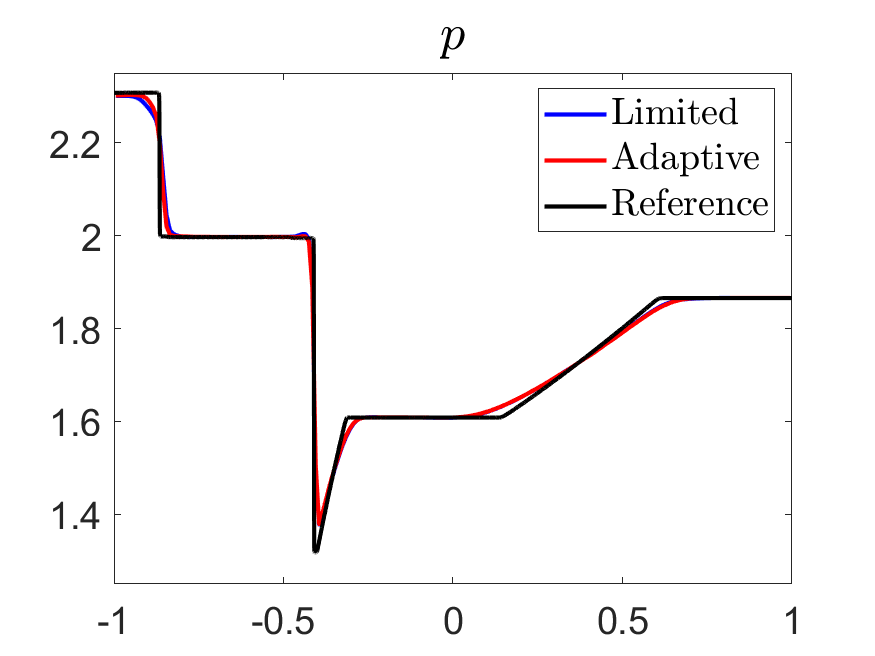}}
\caption{\sf Example 2: Density $\rho$ (left) and pressure $p$ (right) computed by the limited and adaptive schemes.\label{fig8}}
\end{figure}

\begin{figure}[ht!]
\centerline{\includegraphics[trim=0.5cm 0.3cm 0.9cm 0.4cm, clip, width=6.0cm]{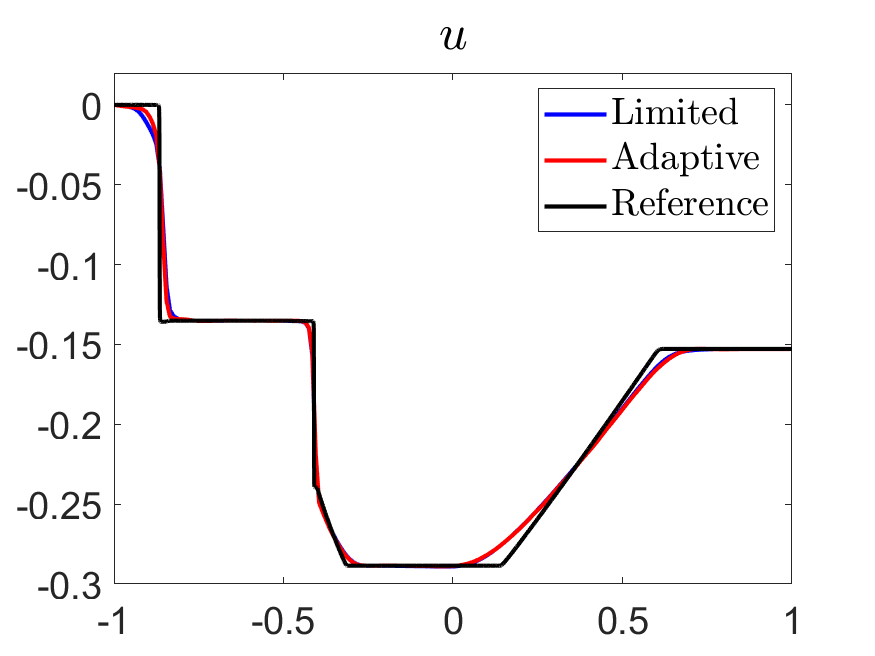}\hspace{1cm}
            \includegraphics[trim=0.5cm 0.3cm 0.9cm 0.4cm, clip, width=6.0cm]{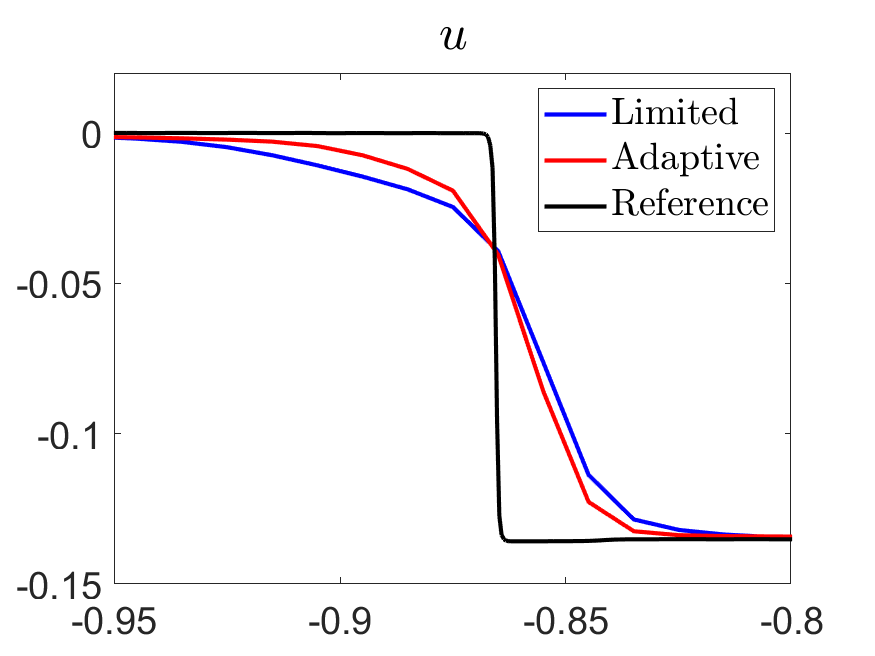}}
\caption{\sf Example 2: Velocity $u$ computed by the limited and adaptive schemes (left) and zoom at $x\in[-0.95,-0.8]$ (right).
\label{fig9}}
\end{figure}

\begin{figure}[ht!]
\centerline{\includegraphics[trim=1.1cm 0.3cm 1.3cm 0.1cm, clip, width=5.5cm]{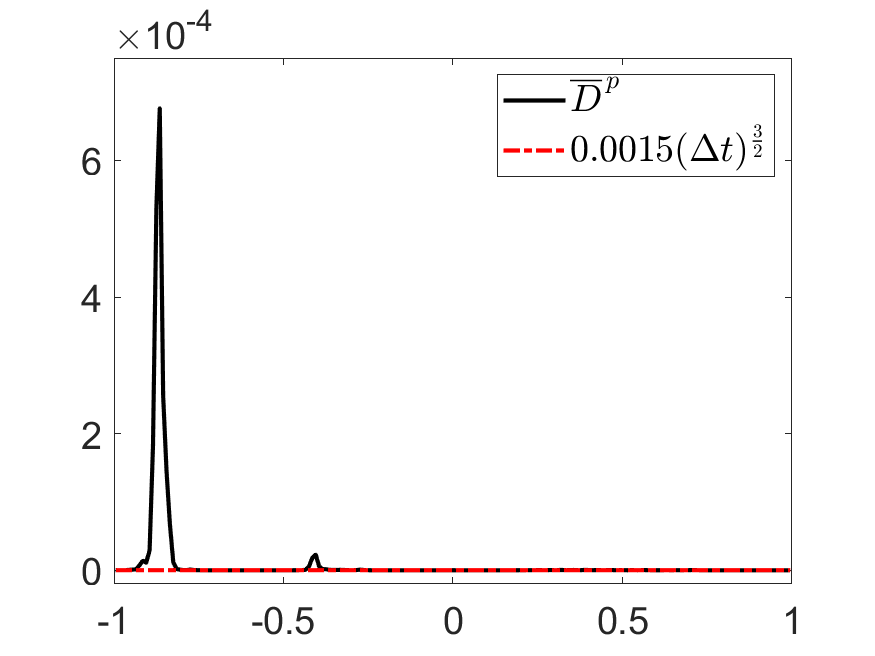}\hspace{1.5cm}
            \includegraphics[trim=1.1cm 0.3cm 1.3cm 0.1cm, clip, width=5.5cm]{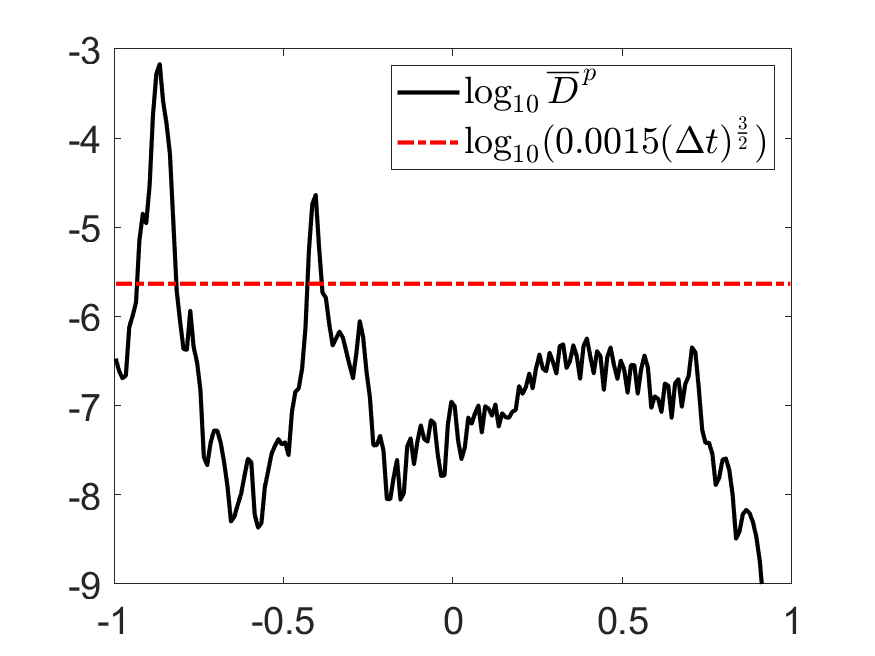}}
\caption{\sf Example 2: $\xbar D^{\,p}$ and $0.0015(\dt)^\frac{3}{2}$ (left) and the corresponding logarithmic quantities (right).
\label{fig10}}
\end{figure}

We compute the numerical solution in the computational domain $[-1,1]$ on the uniform mesh with $\dx=1/100$ and impose the solid wall
boundary conditions at $x=-1$ and free boundary conditions at $x=1$. In Figures \ref{fig8} and \ref{fig9}, we plot the numerical solutions
at the final time $t=3$ obtained by the limited and adaptive (with $\texttt{C}=0.0015$) schemes. These solutions are compared with the
corresponding reference solutions computed by the limited scheme on a much finer mesh with $\dx=1/2000$. In Figure \ref{fig10}, the graphs
of LSI $\xbar D^{\,p}$ and $0.0015(\dt)^\frac{3}{2}$ are depicted along with their logarithm forms. As one can observe, the LSI accurately
captures the position of the shock waves, and the results obtained by the adaptive scheme are a little sharper compared to those obtained by
the limited counterpart.

\paragraph{Example 3---Shock-Entropy Wave Interaction Problem.} In the third example taken from \cite{Shu88}, we consider the shock-entropy
wave interaction problem. The system \eref{3.6}--\eref{3.7} is numerically solved subject to the following initial condition:
\begin{equation*}
(\rho, u,p)(x,0)=\begin{cases}
(1.51695,0.523346,1.805),&x<-4.5,\\
(1+0.1\sin(20x),0,1),&x>-4.5,
\end{cases}
\end{equation*}
which corresponds to a forward-facing shock wave of Mach number 1.1 interacting with high-frequency density perturbations, that is, as the
shock wave moves, the perturbations spread ahead.

\begin{figure}[ht!]
\centerline{\includegraphics[trim=0.8cm 0.3cm 1.2cm 0.4cm, clip, width=6.0cm]{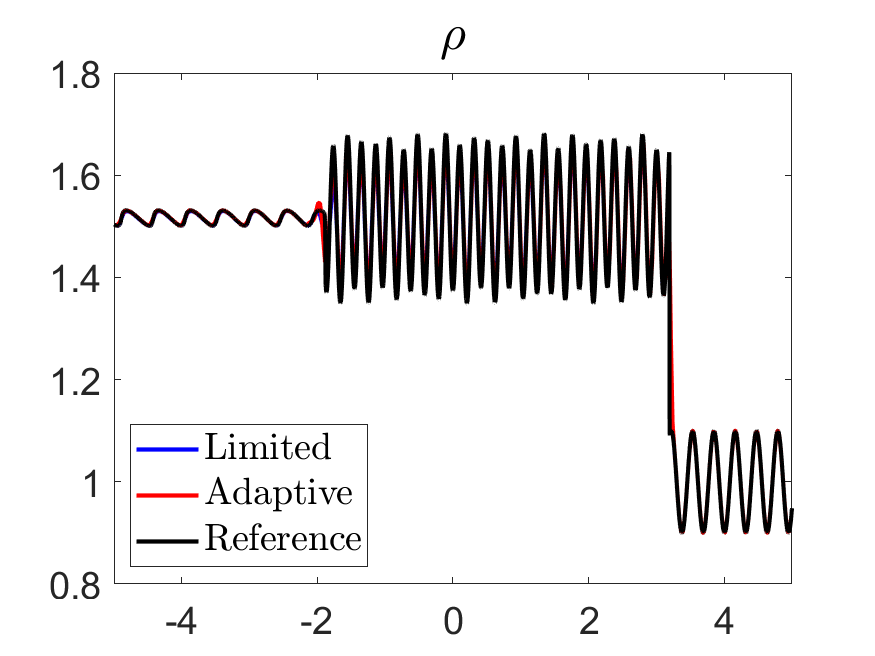}\hspace{1cm}
            \includegraphics[trim=0.8cm 0.3cm 1.2cm 0.4cm, clip, width=6.0cm]{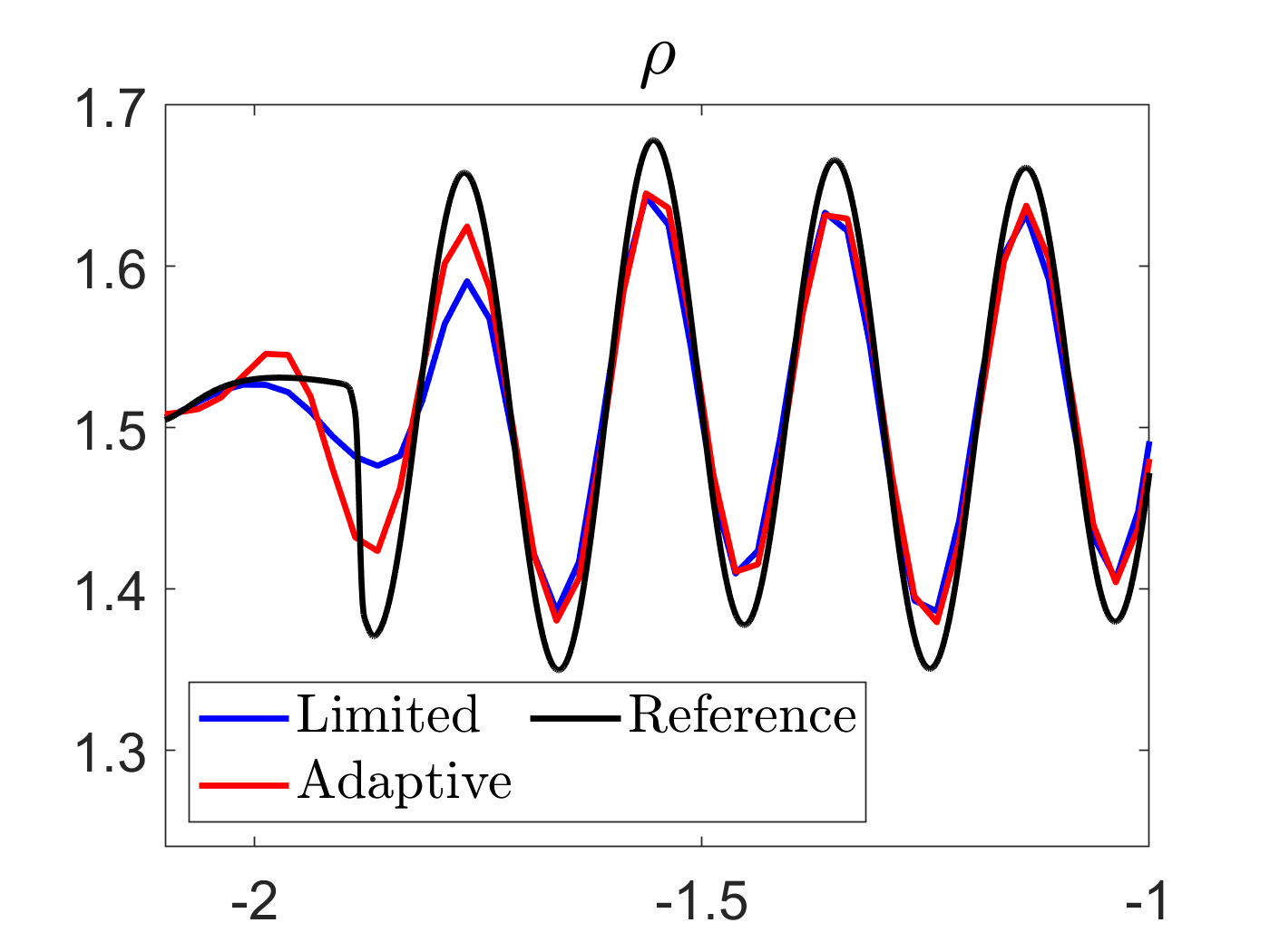}}
\caption{\sf Example 3: Density $\rho$ computed by the limited and adaptive schemes (left) and zoom at $x\in[-2.1,-1]$ (right).\label{fig4}}
\end{figure}

\begin{figure}[ht!]
\centerline{\includegraphics[trim=0.8cm 0.3cm 1.2cm 0.1cm, clip, width=5.5cm]{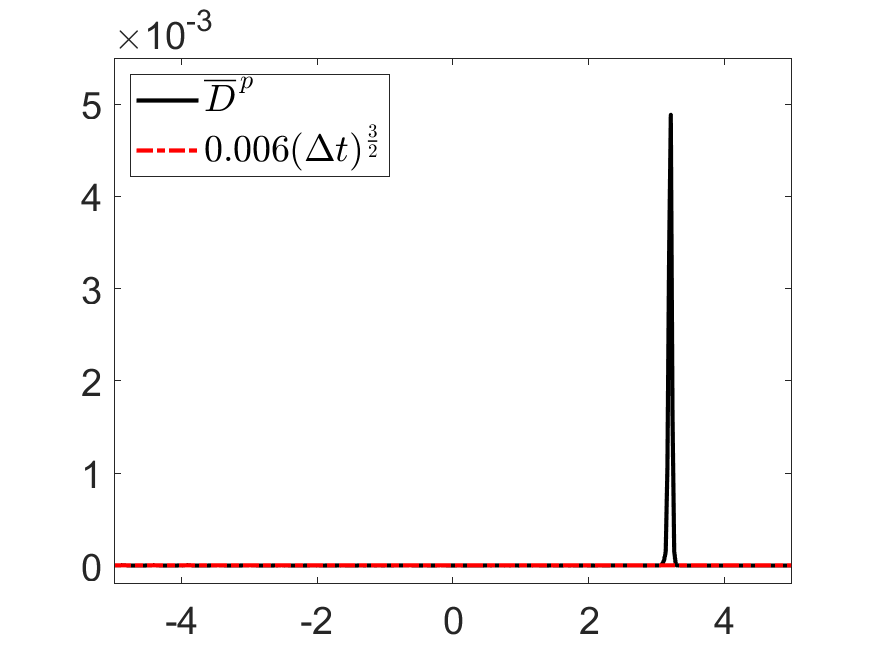}\hspace{1.5cm}
            \includegraphics[trim=0.8cm 0.3cm 1.2cm 0.1cm, clip, width=5.5cm]{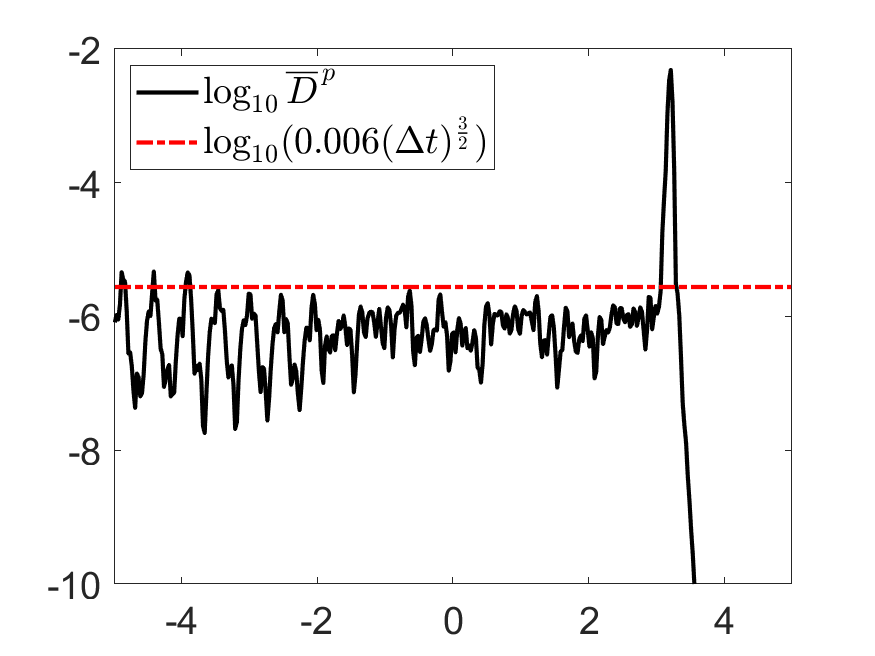}}
\caption{\sf Example 3: $\xbar D^{\,p}$ and $0.006(\dt)^\frac{3}{2}$ (left) and the corresponding logarithmic quantities (right).
\label{fig5}}
\end{figure}

We compute the numerical solution using both the limited and adaptive schemes with $\texttt{C}=0.006$ in the computational domain $[-5,5]$
covered by a uniform mesh with $\dx=1/40$ and implement free boundary conditions. The numerical results at time $t=5$ are presented in
Figure \ref{fig4} along with the corresponding reference solution computed by the limited scheme on a much finer mesh with $\dx=1/800$. As
in the previous example, we also plot (in Figure \ref{fig5}) the graphs of the LSI $\xbar D^{\,p}$ and $0.006(\dt)^\frac{3}{2}$ together
with $\log_{10}\xbar D^{\,p}$ and $\log_{10}(0.006(\dt)^\frac{3}{2})$. One can observe that the LSI can capture the position of the shock
waves accurately, and the results obtained by the adaptive scheme are non-oscillatory and slightly sharper than those obtained by the limited
scheme.

\paragraph{Example 4---Shock-Density Wave Interaction Problem.} In the last 1-D example taken from \cite{Shu89}, we consider the
shock-density wave interaction problem. The initial data,
\begin{equation*}
(\rho,u,p)(x,0)=\begin{cases}
\Big(\dfrac{27}{7},\dfrac{4\sqrt{35}}{9},\dfrac{31}{3}\Big),&x<-4,\\[0.8ex]
(1+0.2\sin(5x),0,1),&x>-4,
\end{cases}
\end{equation*}
are prescribed in the computational domain $[-5,15]$ subject to the free boundary conditions.

\begin{figure}[ht!]
\centerline{\includegraphics[trim=0.8cm 0.3cm 0.8cm 0.4cm, clip, width=6.0cm]{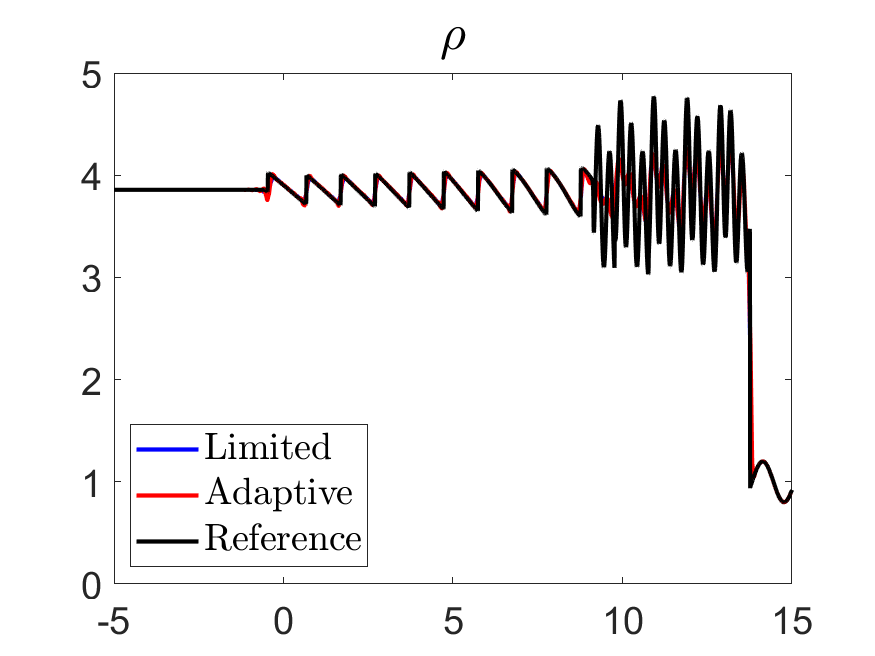}\hspace{1cm}
            \includegraphics[trim=0.8cm 0.3cm 0.8cm 0.4cm, clip, width=6.0cm]{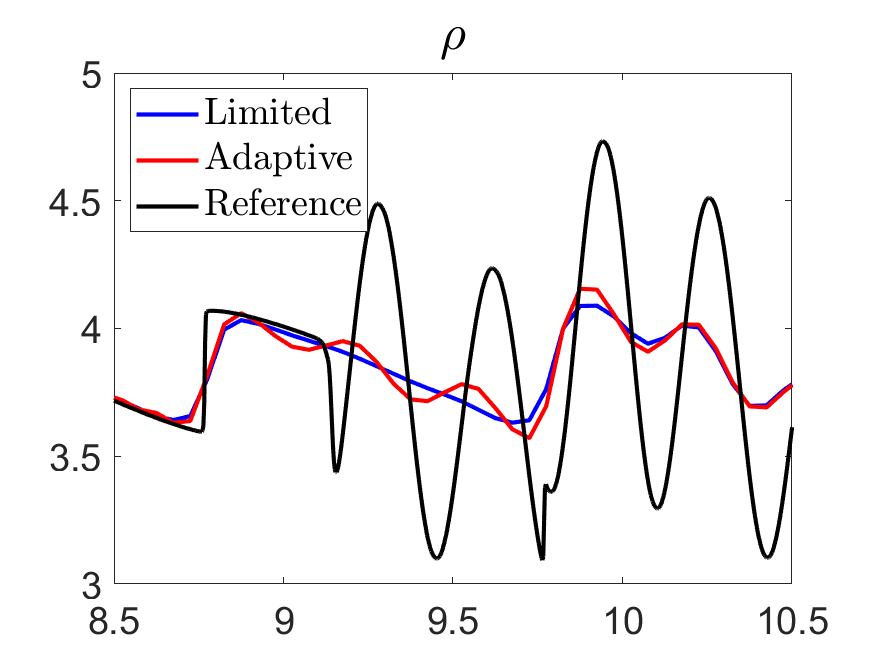}}
\caption{\sf Example 4: Density $\rho$ computed by the limited and adaptive schemes (left) and zoom at $x\in[8,10]$ (right).\label{fig9a}}
\end{figure}

\begin{figure}[ht!]
\centerline{\includegraphics[trim=0.1cm 0.3cm 1.1cm 0.6cm, clip, width=6.0cm]{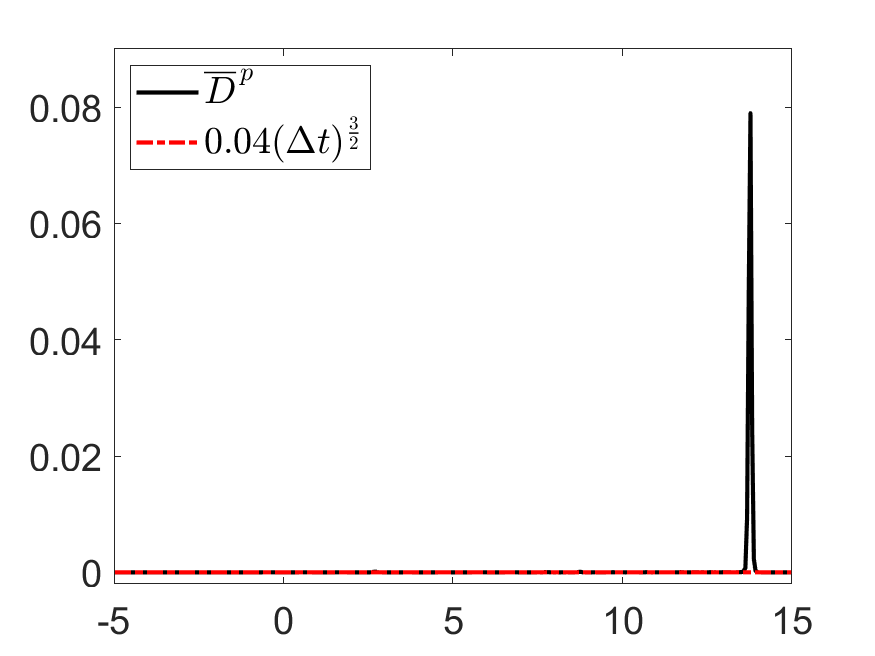}\hspace{1cm}
            \includegraphics[trim=0.1cm 0.3cm 1.1cm 0.6cm, clip, width=6.0cm]{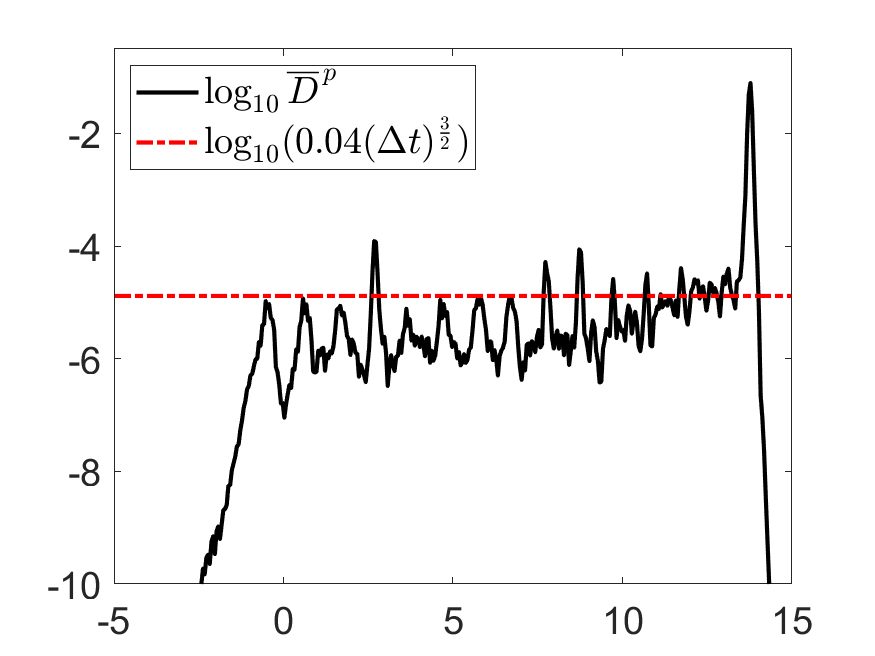}}
\caption{\sf Example 4: $\xbar D^{\,p}$ and $0.04(\dt)^\frac{3}{2}$ (left) and the corresponding logarithmic quantities (right).
\label{fig9b}}
\end{figure}

We compute the numerical solution by the limited and adaptive (with $\texttt{C}=0.04$) schemes on the uniform mesh with $\dx=1/20$ until the
final time $t=5$ and present the obtained numerical results in Figure \ref{fig9a} together with the corresponding reference computed by the
limited scheme on a much finer mesh with $\dx=1/400$. We also plot $\xbar D^{\,p}$ and $0.04(\dt)^\frac{3}{2}$ together with
$\log_{10}\xbar D^{\,p}$ and $\log_{10}(0.04(\dt)^\frac{3}{2})$ in Figure \ref{fig9b}. It can be seen clearly that the LSI can accurately
capture the position of the shock waves, and the adaptive scheme produces slightly sharper results compared to those obtained by the limited
scheme.

\begin{rmk}
It is instructive to compare the computational costs of the studied limited and adaptive A-WENO schemes. To this end, we have measured the
CPU times consumed by both schemes. The results obtained for the four studied 1-D examples are reported in Table \ref{tab2}, where we show
the relative CPU time consumption of the adaptive A-WENO scheme relative to the fully limited one. As one can see, the proposed adaptive
the scheme is more efficient than the fully limited one. Notice that the numbers in Table \ref{tab2} are different as the part of the
computational domain indicated as ``rough'' varies. The CPU times for the adaptive scheme also depend on the values of $\texttt{C}$: the use
of larger $\texttt{C}$ leads to a more efficient but potentially more oscillatory adaptive A-WENO scheme.
\begin{table}[ht!]
\centering
\begin{tabular}{|c|c|c|c|}
\hline
Example 1&Example 2&Example 3&Example 4\\
\hline
66\% &76\% &66\% &67\%\\
\hline
\end{tabular}
\caption{\sf Examples 1--4: CPU times consumed by the adaptive A-WENO scheme relative to the fully limited A-WENO scheme.\label{tab2}}
\end{table}
\end{rmk}

\subsection{2-D Examples}\label{sec42}
In this section, we demonstrate the performance of the proposed adaptive A-WENO scheme on several examples for the 2-D Euler equations of
gas dynamics, which read as
\begin{equation}
\begin{aligned}
&\rho_t+(\rho u)_x+(\rho v)_y=0,\\
&(\rho u)_t+(\rho u^2 +p)_x+(\rho uv)_y=0,\\
&(\rho v)_t+(\rho uv)_x+(\rho v^2+p)_y=0,\\
&E_t+\left[u(E+p)\right]_x+\left[v(E+p)\right]_y=0,
\end{aligned}
\label{4.3}
\end{equation}
where $v$ is the $y$-component of the velocity, and the rest of the notations are the same as in the 1-D case. The system is completed through
the following equations of state:
\begin{equation}
p=(\gamma-1)\Big[E-\frac{\rho}{2}(u^2+v^2)\Big].
\label{4.4}
\end{equation}

\paragraph{Example 5---2-D Riemann Problem.} In the first 2-D example, we consider Configuration 3 of the 2-D Riemann problems from
\cite{Kurganov02} (see also \cite{Schulz93,Schulz93a,Zheng01}) with the following initial conditions:
\begin{equation*}
(\rho,u,v,p)(x,y,0)=\begin{cases}
(1.5,0,0,1.5),&x>1,~y>1,\\
(0.5323,1.206,0,0.3),&x<1,~y>1,\\
(0.138,1.206,1.206,0.029),&x<1,~y<1,\\
(0.5323,0,1.206,0.3),&x>1,~y<1.
\end{cases}
\end{equation*}

\begin{figure}[ht!]
\centerline{\includegraphics[trim=5.0cm 3.4cm 1.8cm 2.4cm, clip, width=14cm]{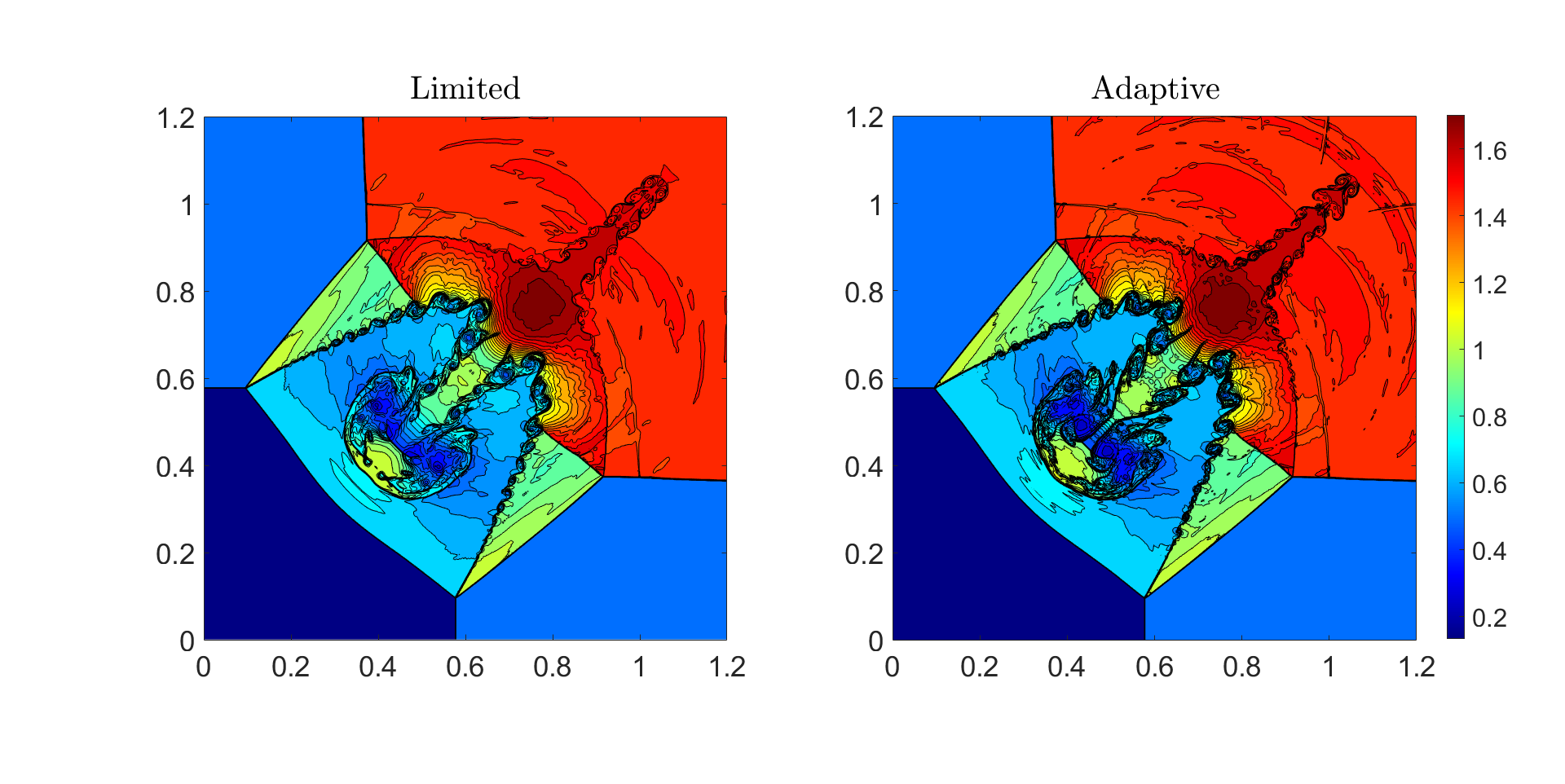}}
\caption{\sf Example 5: Density $\rho$ computed by the limited (left) and adaptive (right) schemes.\label{fig11a}}
\end{figure}

We compute the numerical solution until the final time $t=1$ by the limited and adaptive (with $\texttt{C}=3$) schemes on the uniform mesh
with $\dx=\dy=3/2500$ in the computational domain $[0,1.2]\times[0,1.2]$ subject to the free boundary conditions. The obtained results are
presented in Figure \ref{fig11a}, where one can see that the adaptive scheme outperforms the limited one as it better captures
the sideband instability of the jet in the zones of strong along-jet velocity shear and the instability along the jet’s neck.

\begin{figure}[ht!]
\centerline{\includegraphics[trim=2.0cm 0.3cm 2.1cm 0.6cm, clip, width=5.5cm]{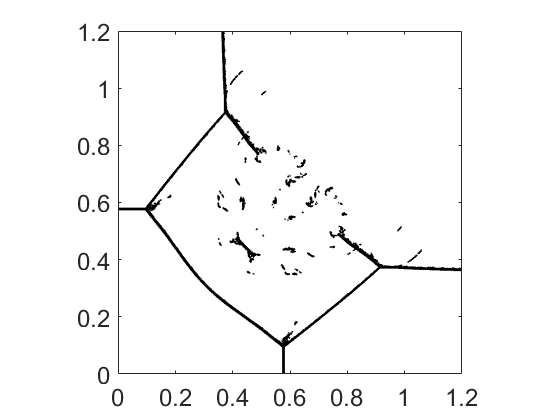}}
\caption{\sf Example 5: The limited WENO-Z interpolation is used only in the part of the computational domain indicated by the black color.
\label{fig10a}}
\end{figure}

In Figure \ref{fig10a}, we show the regions which the LSI detected as ``rough'' at the final time. As one can see, the limited
WENO-Z interpolation is used only in a small part of the computational domain, mostly around the shocks.

\paragraph{Example 6---Explosion Problem.} In this example, we consider the explosion problem taken from \cite{Liska03,Tor}. This is a
circularly symmetric problem with the following initial conditions:
\begin{equation}
(\rho,u,v,p)(x,y,0)=\begin{cases}
(1,0,0,1),&x^2+y^2<0.16,\\
(0.125,0,0,0.1),&\mbox{otherwise}.
\end{cases}
\label{5.3}
\end{equation}
We numerically solve the initial value problem \eref{4.3}--\eref{5.3} in the first quadrant, more precisely in the computational domain
$[0,1.5]\times[0,1.5]$ with the solid wall boundary conditions imposed at $x=0$ and $y=0$ and the free boundary conditions set at $x=1.5$
and $y=1.5$.

\begin{figure}[ht!]
\centerline{\includegraphics[trim=5.0cm 3.5cm 1.5cm 2.5cm, clip, width=14cm]{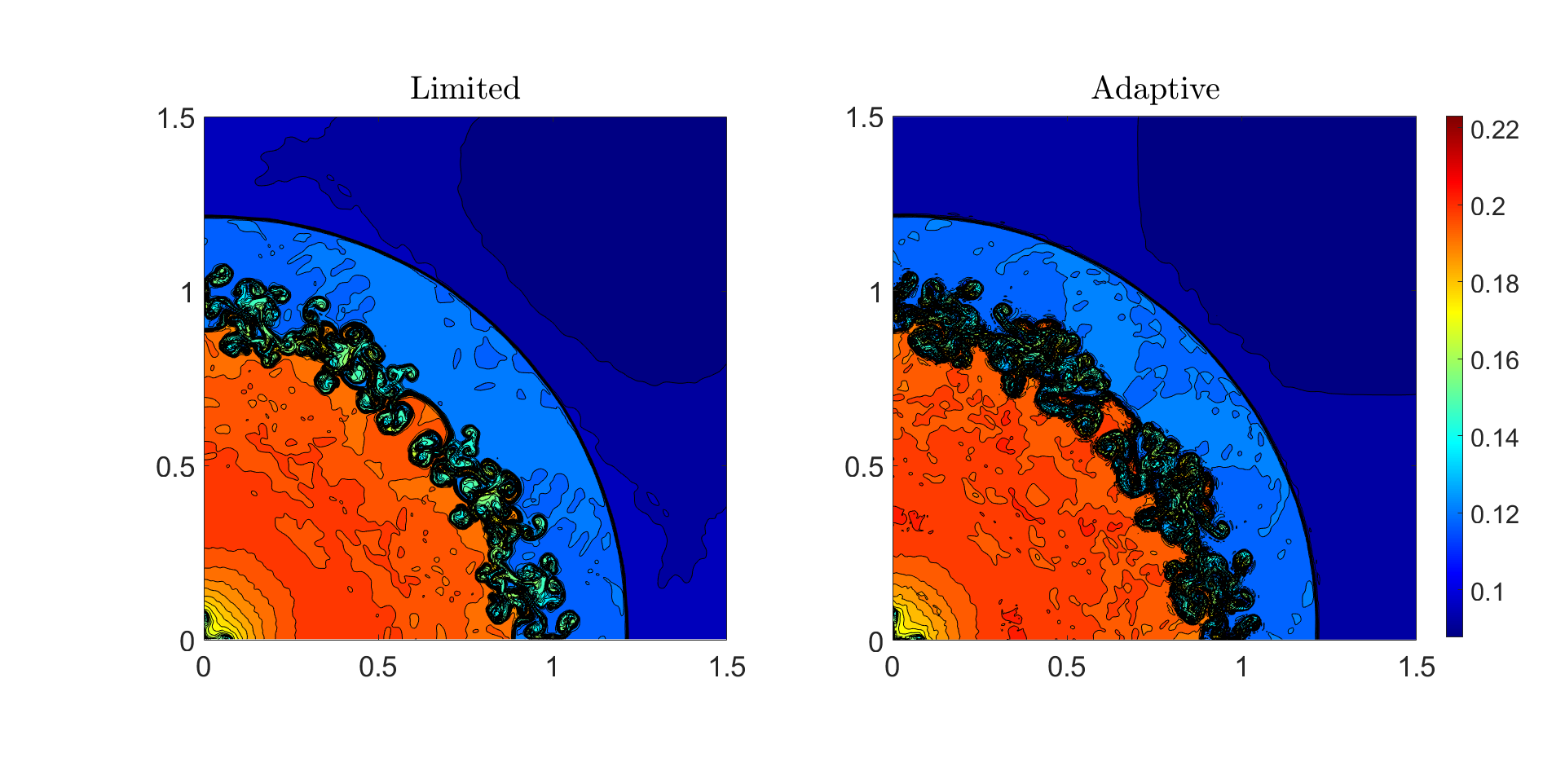}}
\caption{\sf Example 6: Density $\rho$ computed by the limited (left) and adaptive (right) schemes.\label{fig11}}
\end{figure}

In Figure \ref{fig11}, the numerical solutions computed by the limited and adaptive (with $\texttt{C}=1$) schemes on the uniform mesh with
$\dx=\dy=3/800$ are plotted at the final time $t=3.2$. The presented results clearly illustrate the advantage of the adaptive approach over
the fully limited one, as the contact curve captured by the adaptive scheme is much ``curlier'' and the mixing layer is much ``wider''.

\begin{figure}[ht!]
\centerline{\includegraphics[trim=5.2cm 0.9cm 5.5cm 1.6cm, clip, width=5.5cm]{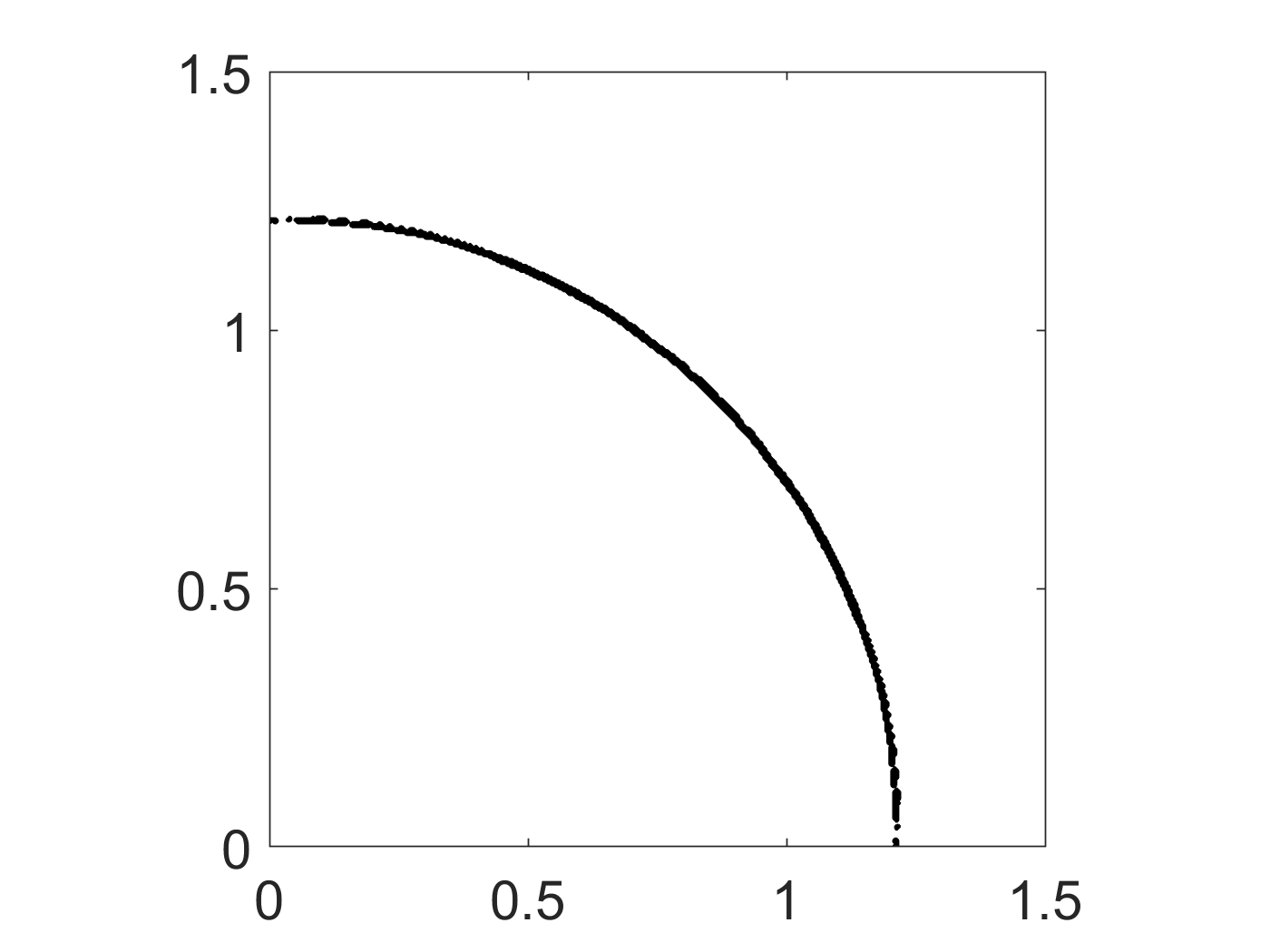}}
\caption{\sf Example 6: The limited WENO-Z interpolation is used only in the part of the computational domain indicated by the black
color.\label{fig11aa}}
\end{figure}

In Figure \ref{fig11aa}, we show the regions detected at the final time by the LSI as ``rough'' and demonstrate that in this example, the
limiting is only used along the circular shock.

\paragraph{Example 7---Implosion Problem.}
In this example taken from \cite{Liska03}, we consider the implosion problem with the following initial conditions:
\begin{equation}
(\rho,u,v,p)(x,y,0)=\begin{cases}
(0.125,0,0,0.14),&|x|+|y|<0.15,\\
(1,0,0,1),&\mbox{otherwise}.
\end{cases}
\label{5.4a}
\end{equation}
prescribed in $[-0.3,0.3]\times[-0.3,0.3]$ subject to the solid wall boundary conditions. Due to the symmetry, we numerically solve the
initial-boundary value problem \eref{4.3}, \eref{4.4} and \eref{5.4a} in the first quadrant only, more precisely in the computational
domain $[0,0.3]\times[0,0.3]$ and impose the solid wall boundary conditions at $x=0$ and $y=0$.

\begin{figure}[ht!]
\centerline{\includegraphics[trim=5.0cm 3.4cm 1.9cm 2.5cm, clip, width=14cm]{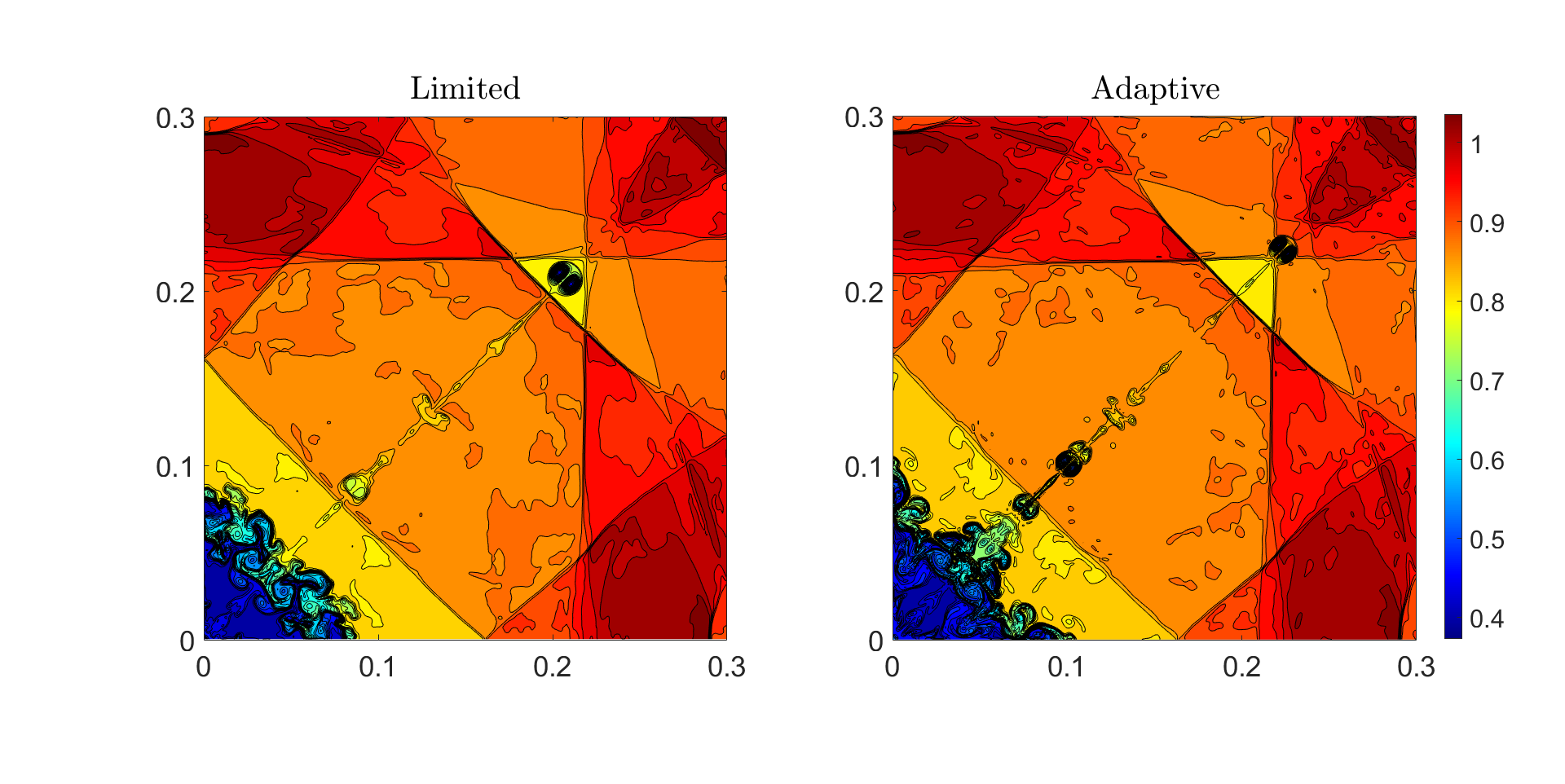}}
\caption{\sf Example 7: Density $\rho$ computed by the limited (left) and adaptive (right) schemes.\label{fig12}}
\end{figure}

In Figure \ref{fig12}, the numerical solutions computed by the limited and adaptive (with $\texttt{C}=3$) schemes on the uniform mesh with
$\dx=\dy=3/4000$ are plotted at the final time $t=2.5$. As one can observe, the jet generated by the adaptive scheme propagates further in
the direction of $y=x$ than the jet produced by the limited scheme, clearly indicating that the adaptive scheme is substantially less
dissipative than the limited scheme.

\begin{figure}[ht!]
\centerline{\includegraphics[trim=5.1cm 0.9cm 5.5cm 1.6cm, clip, width=5.5cm]{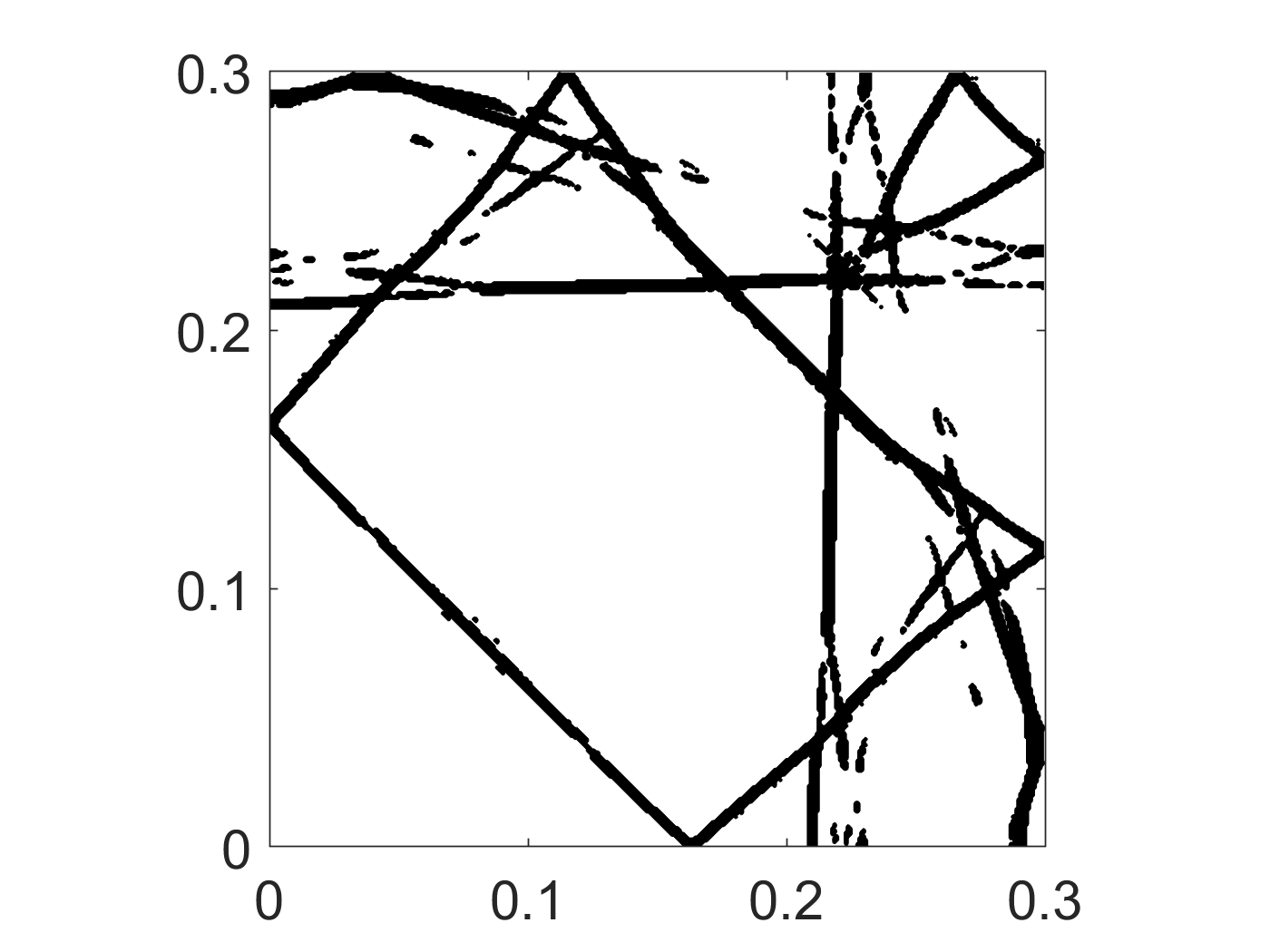}}
\caption{\sf Example 7: The limited WENO-Z interpolation is used only in the part of the computational domain indicated by the black color.
\label{fig12aa}}
\end{figure}

The domain where the limiters have been used at the final time is presented in Figure \ref{fig12aa}, where one can see how the proposed
LSI identifies ``rough'' areas.

\begin{rmk}
It is easy to show that the solution of the studied initial-boundary value problem is symmetric with respect to the axis $y=x$. It is
well-known, however, that this symmetry may be destroyed by the roundoff errors when the solution is computed by a low-dissipative
high-order scheme. In order to prevent the loss of symmetry, we have used a very simple strategy introduced in \cite{WDGK2020}: upon
completion of each time evolution step, we replace the computed point values $\bm U_{j,k}$ with $\widehat{\bm U}_{j,k}$, where
\begin{equation*}
\widehat\rho_{j,k}:=\frac{\rho_{j,k}+\rho_{k,j}}{2},~~(\widehat{\rho u})_{j,k}:=\frac{(\rho u)_{j,k}+(\rho v)_{k,j}}{2},~~
(\widehat{\rho v})_{j,k}:=\frac{(\rho v)_{j,k}+(\rho u)_{k,j}}{2},~~\widehat E_{j,k}:=\frac{E_{j,k}+E_{k,j}}{2},
\end{equation*}
for all $j,k$. For more sophisticated symmetry enforcement techniques, we refer the reader to, e.g., \cite{DLGW,Fleischman19,DLWG,WTX}.
\end{rmk}

\paragraph{Example 8---KH Instability.} In this example taken from \cite{Fjordholm16,Panuelos20}, we study the KH instability with the
following initial conditions:
\begin{equation*}
(\rho(x,y,0), u(x,y,0))=\begin{cases}
(1,-0.5+0.5e^{(y+0.25)/L}),&y\in[-0.5,-0.25),\\
(2,0.5-0.5e^{(-y-0.25)/L}),&y\in[-0.25,0),\\
(2,0.5-0.5e^{(y-0.25)/L}),&y\in[0,0.25),\\
(1,-0.5+0.5e^{(-y+0.25)/L}),&y\in[0.25,0.5),
\end{cases}
\end{equation*}
\begin{equation*}
v(x,y,0)=0.01\sin(4\pi x),\quad p(x,y,0)\equiv1.5,
\end{equation*}
where $L$ is a smoothing parameter (here, we take $L=0.00625$) corresponding to a thin shear interface with a perturbed vertical velocity
field $v$ in the conducted simulations. We impose the 1-periodic boundary conditions in both the $x$- and $y$-directions, and take the
computational domain to be $[-0.5,0.5]\times[-0.5,0.5]$.

\begin{figure}[ht!]
\centerline{\includegraphics[trim=5.2cm 6.6cm 2.2cm 5.4cm, clip, width=17cm]{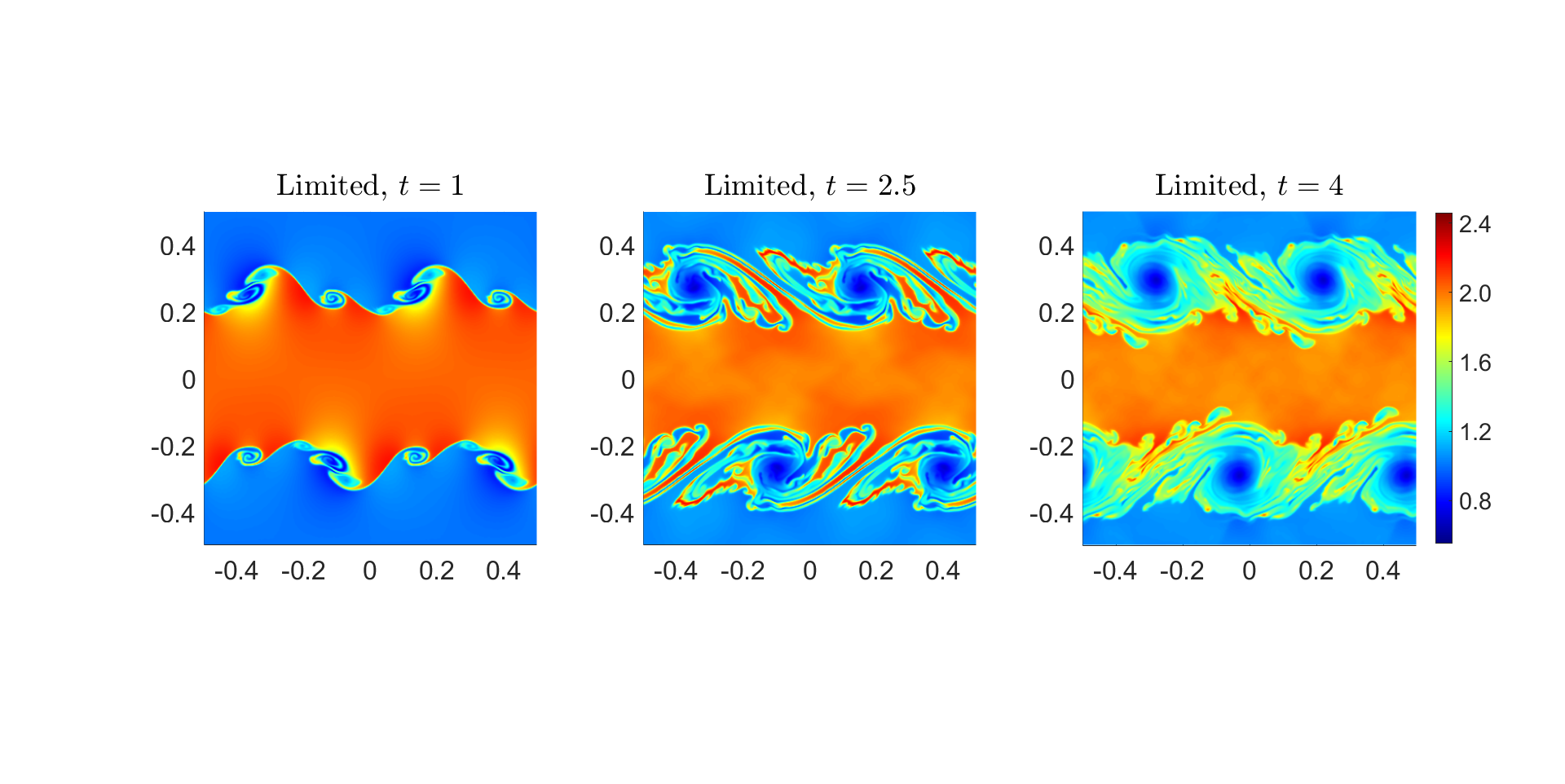}}
\vskip8pt
\centerline{\includegraphics[trim=5.2cm 6.6cm 2.2cm 5.4cm, clip, width=17cm]{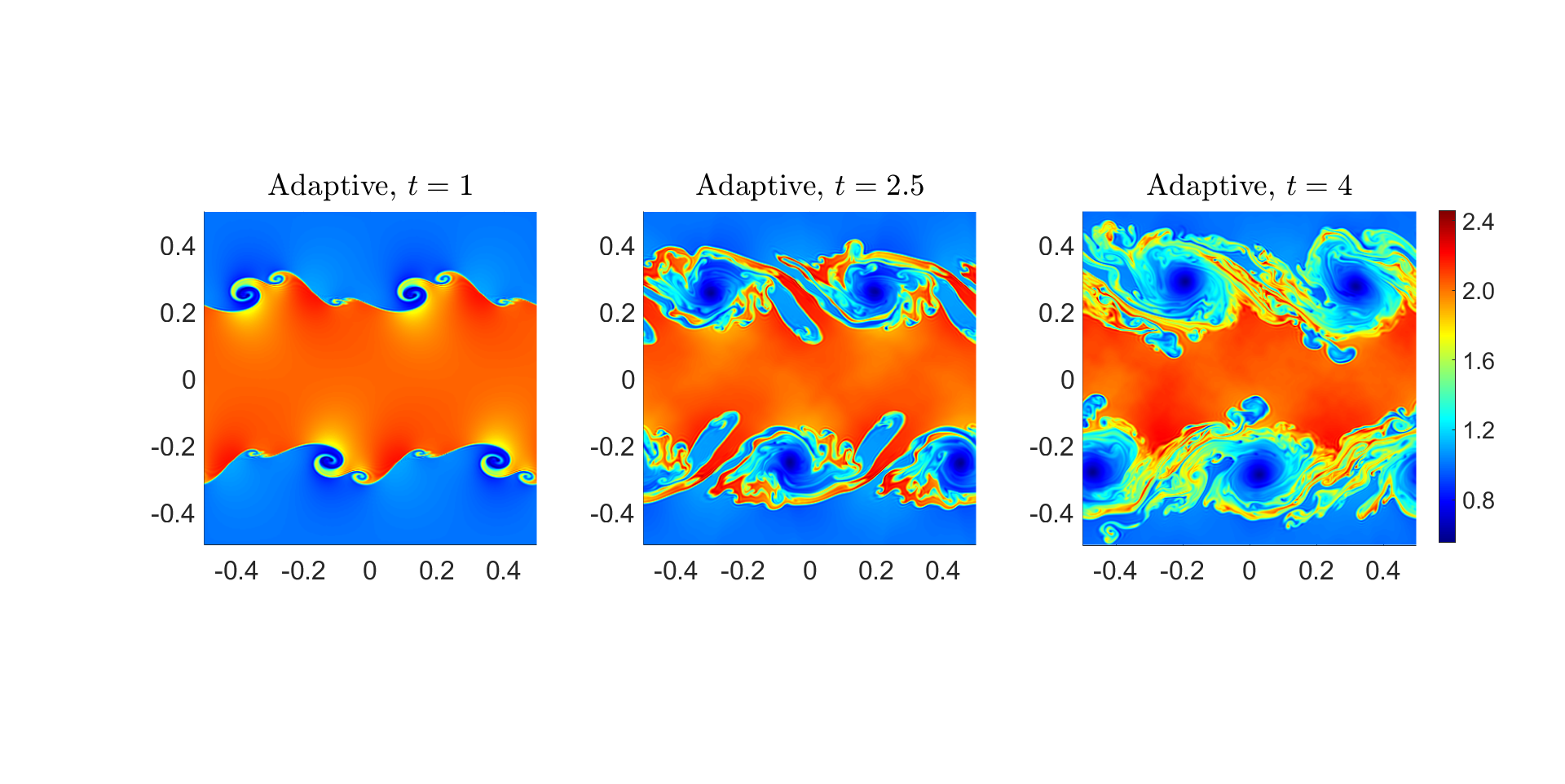}}
\caption{\sf Example 8: Time snapshots of the density $\rho$ computed by the limited (top row) and adaptive (bottom row) schemes at $t=1$
(left column), 2.5 (middle column), and 4 (right column).\label{fig13}}
\end{figure}

We compute the numerical solution until the final time $t=4$ by the limited and adaptive (with $\texttt{C}=1$) schemes on the uniform mesh
with $\dx=\dy=1/400$. The numerical results at $t=1$, 2.5, and 4 are presented in Figure \ref{fig13}. As one can see, at the early time
$t=1$, the vortex sheets in the limited and adaptive results are quite different, and it is hard to draw a definite conclusion based on these
results. However, at later times $t=2.5$ and 4, the adaptive scheme produces more complicated vortices and turbulent mixing, which
indicates that the adaptive scheme contains less numerical dissipation than the limited scheme.

\begin{figure}[ht!]
\hspace*{-0.5cm}\centerline{\includegraphics[trim=9.4cm 12.5cm 9.6cm 11.0cm, clip, width=16cm]{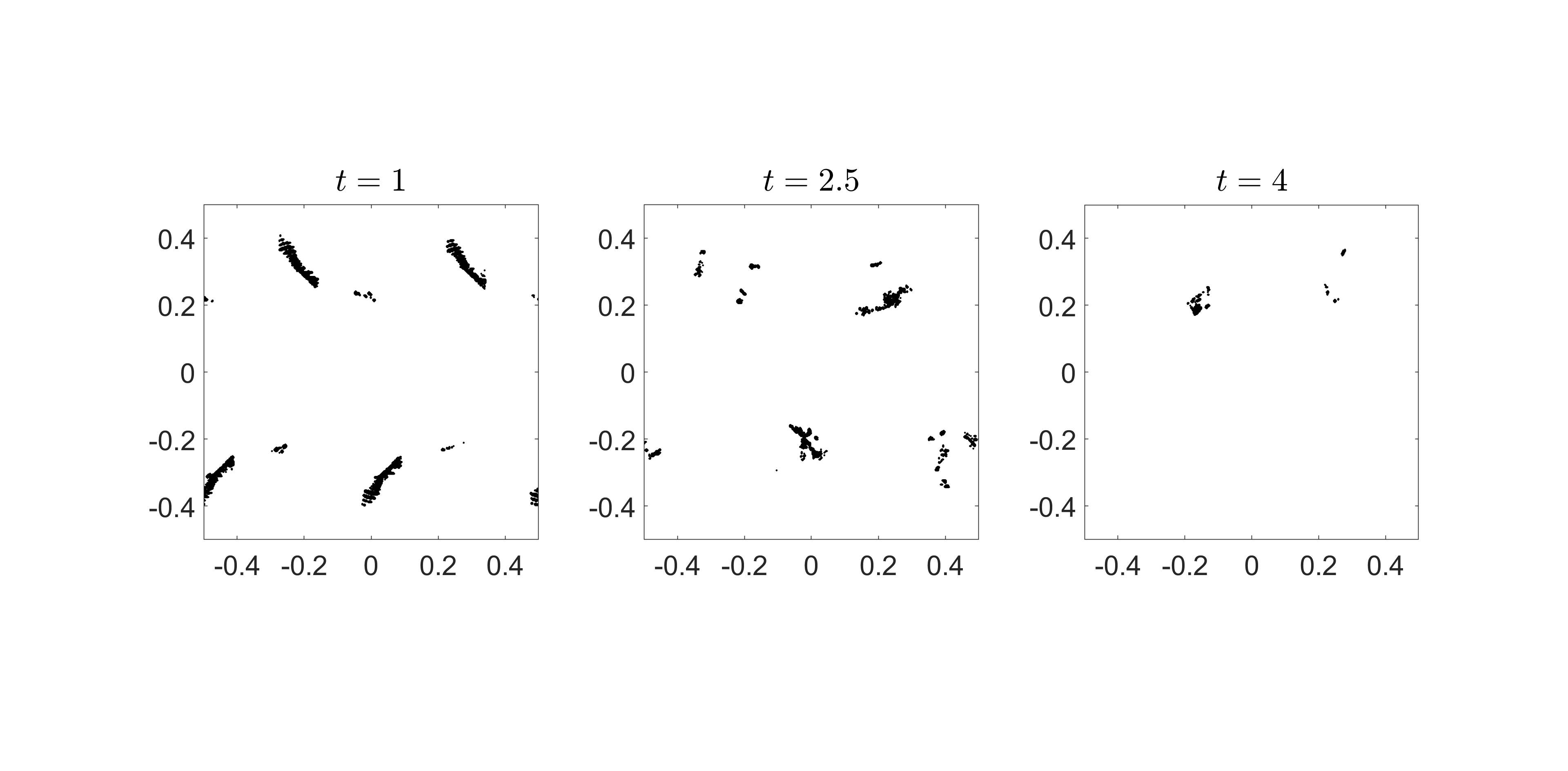}}
\caption{\sf Example 8: The limited WENO-Z interpolation is used only in the part of the computational domain indicated by the black color.
\label{fig13aa}}
\end{figure}

In addition, in Figure \ref{fig13aa}, we plot the solution regions to show that the limiters have been used in a very small part of the
computational domain, especially at $t=4$.

\begin{figure}[ht!]
\centerline{\includegraphics[trim=5.2cm 6.6cm 2.4cm 5.6cm, clip, width=17cm]{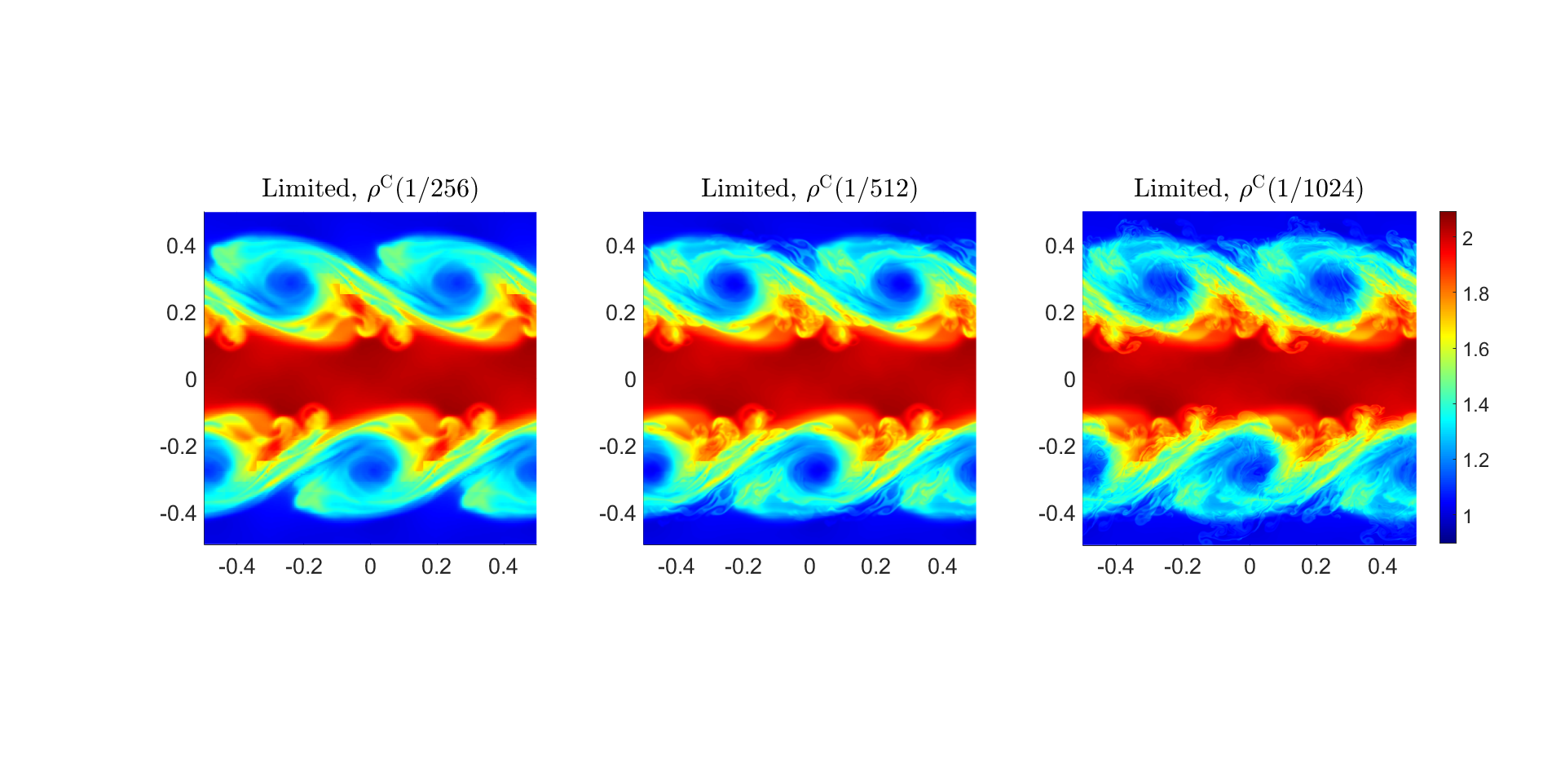}}
\vskip8pt
\centerline{\includegraphics[trim=5.2cm 6.6cm 2.4cm 5.6cm, clip, width=17cm]{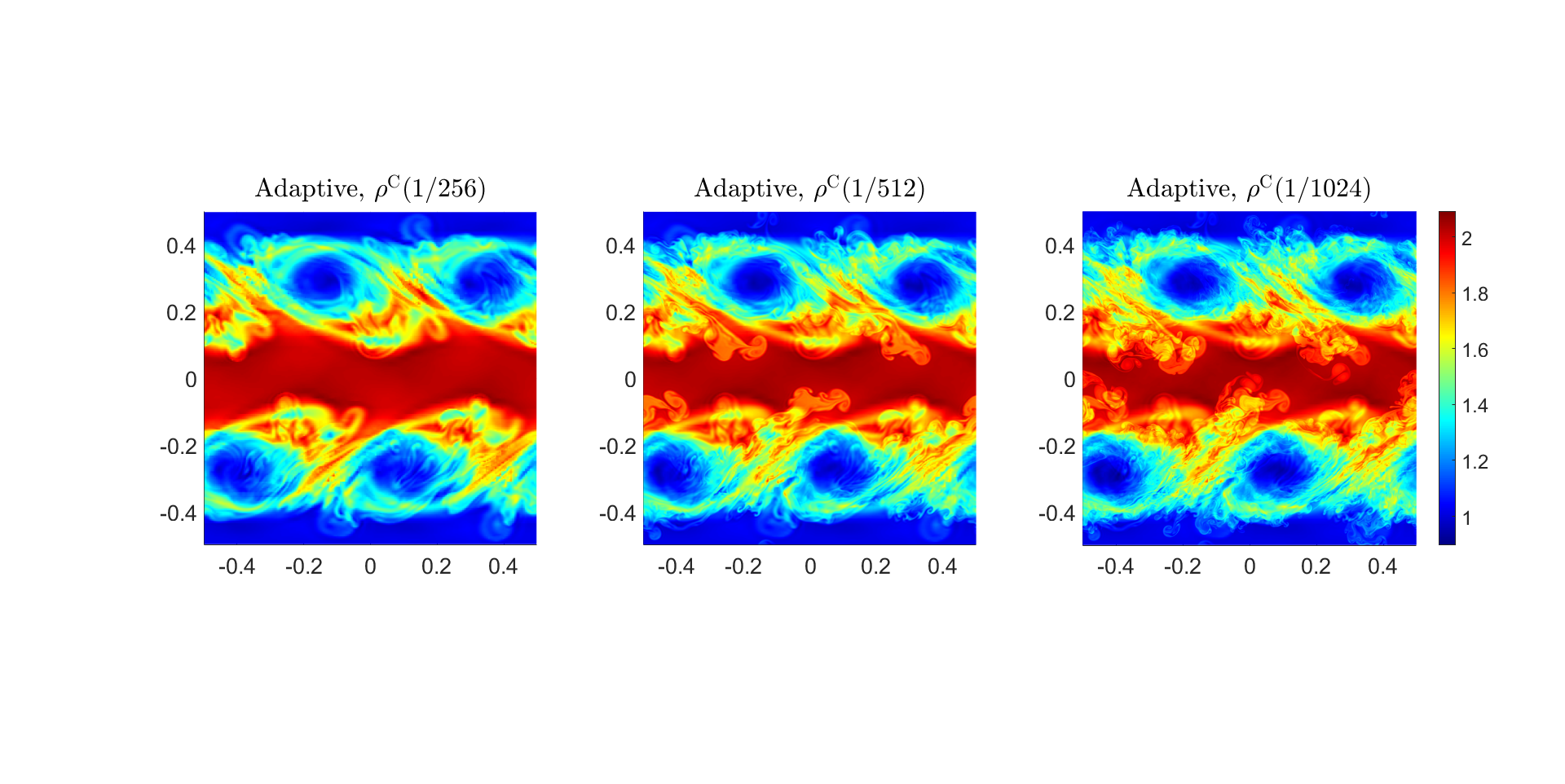}}
\caption{\sf Example 8: Ces\`aro averages of the density $\rho^{\rm C}(1/2^m)$ computed by the limited (top row) and adaptive
(bottom row) schemes for $m=8$ (left column), 9 (middle column), and 10 (right column).\label{fig13bb}}
\end{figure}

It should also be noted that, as is known, the numerical solutions of the KH instability problem do not converge in the strong sense when
the mesh is refined. In fact, the limiting solution is not a weak solution but a dissipative weak solution; see \cite{Lukacova_book} for
more details. Thus, to approximate the limiting solution, we compute the Ces\`aro averages of the densities obtained at the final
time $t=4$ by the limited and adaptive schemes. To this end, we first introduce a sequence of meshes with the cells of size $1/2^n$,
$n=5,\ldots,10$, and denote by $\rho(1/2^n)$ the density computed on the corresponding mesh. We then project the obtained coarser mesh
solutions with $n=5,\ldots,m-1$ onto the finer mesh with $n=m$ (the projection is carried out using the dimension-by-dimension
WENO-Z interpolation of the density field) and denote the obtained densities still by $\rho(1/2^n)$, $n=5,\ldots,m$. After this, the
Ces\`aro averages are computed by
\begin{equation}
\rho^{\rm C}(1/2^m)=\frac{\rho(1/2^5)+\cdots+\rho(1/2^m)}{m-4},\quad m=8,9,10.
\label{5.4}
\end{equation}
In Figure \ref{fig13bb}, we plot the computed averages at time $t=4$. One can observe the superiority of the results
obtained by the adaptive scheme when it comes to resolving complicated structures.

\paragraph{Example 9---RT Instability}
In the last example, we investigate the RT instability. It is a physical phenomenon occurring when a layer of heavier fluid is placed on
top of a layer of lighter fluid. To this end, we first modify the 2-D Euler equations of gas dynamics \eref{4.3}--\eref{4.4} by
adding the gravitational source terms acting in the positive direction of the $y$-axis into the RHS of the system:
\begin{equation*}
\begin{aligned}
&\rho_t+(\rho u)_x+(\rho v)_y=0,\\
&(\rho u)_t+(\rho u^2 +p)_x+(\rho uv)_y=0,\\
&(\rho v)_t+(\rho uv)_x+(\rho v^2+p)_y=\rho,\\
&E_t+\left[u(E+p)\right]_x+\left[v(E+p)\right]_y=\rho v,
\end{aligned}
\end{equation*}
and then use the setting from \cite{Shi03,Wang20} with the following initial conditions:
\begin{equation*}
(\rho,u,v,p)(x,y,0)=\begin{cases}
(2,0,-0.025c\cos(8\pi x),2y+1),&y<0.5,\\
(1,0,-0.025c\cos(8\pi x),y+1.5),&\mbox{otherwise},
\end{cases}
\end{equation*}
where $c:=\sqrt{\gamma p/\rho}$ is the speed of sound. The solid wall boundary conditions are imposed at $x=0$ and $x=0.25$, and the
following Dirichlet boundary conditions are specified at the top and bottom boundaries:
$$
(\rho,u,v,p)(x,1,t)=(1,0,0,2.5),\quad(\rho,u,v,p)(x,0,t)=(2,0,0,1).
$$

\begin{figure}[ht!]
\centerline{\includegraphics[trim=5.0cm 4.5cm 2.6cm 3.5cm, clip, width=17.5cm]{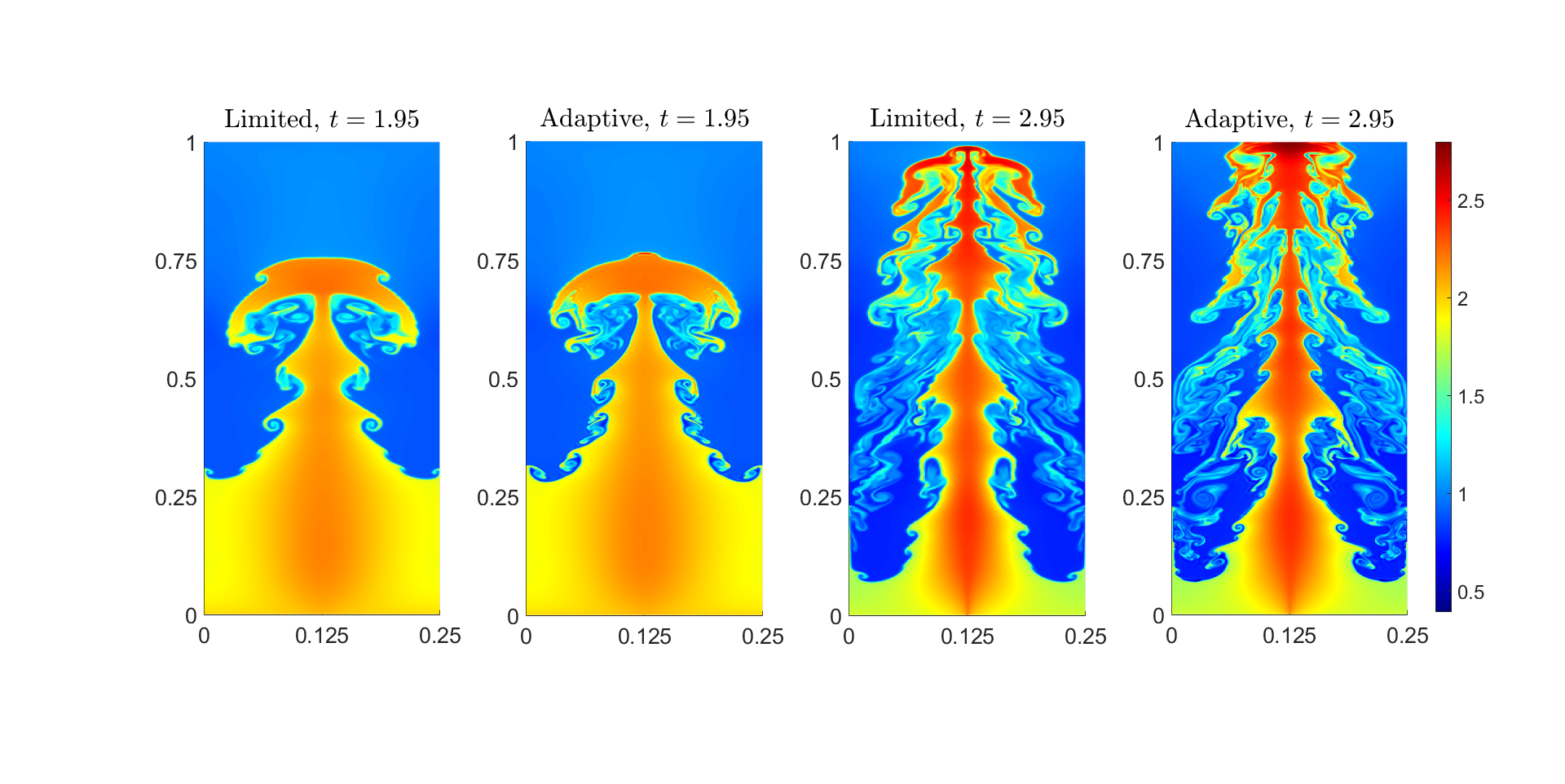}}
\caption{\sf Example 9: Density $\rho$ computed by the limited and adaptive schemes at $t=1.95$ and 2.95.\label{fig14}}
\end{figure}

We compute the numerical solution until the final time $t=2.95$ by the limited and adaptive (with $\texttt{C}=2$) schemes on the
computational domain $[0,0.25]\times[0,1]$ on the uniform mesh with $\dx=\dy=1/800$. The numerical results at $t=1.95$ and 2.95 are
presented in Figure \ref{fig14}. As we can see, there are pronounced differences between the limited and adaptive solutions. Therefore, one
can conclude that the adaptive scheme achieves a much better resolution, which again demonstrates that the adaptive scheme is less
dissipative than the limited scheme.

\begin{figure}[ht!]
\centerline{\includegraphics[trim=4.0cm 0.3cm 4.3cm 0.1cm, clip, width=4.cm]{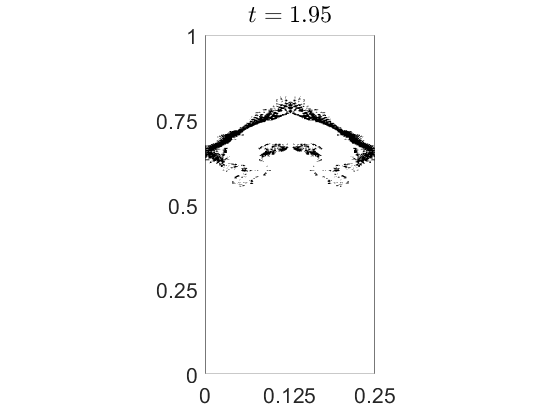}\hspace*{1.0cm}
            \includegraphics[trim=4.0cm 0.3cm 4.3cm 0.1cm, clip, width=4.cm]{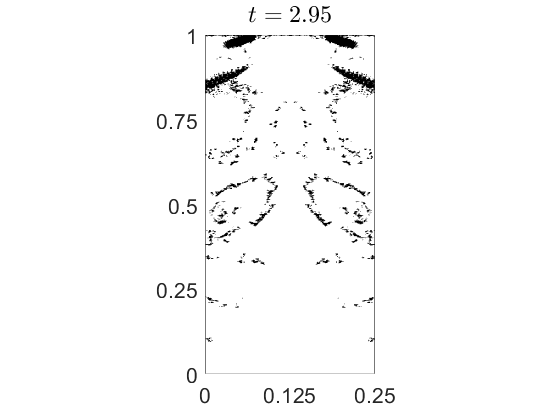}}
\caption{\sf Example 9: The limited WENO-Z interpolation is used only in the part of the computational domain indicated by the black color.
\label{fig14a}}
\end{figure}

In Figure \ref{fig14a}, we show the regions which the LSI detected as ``rough'' at the final time. As one can see, the limited
WENO-Z interpolation is only used in a relatively small part of the computational domain.

\begin{figure}[ht!]
\centerline{\includegraphics[trim=5.1cm 4.5cm 2.5cm 3.5cm, clip, width=17.5cm]{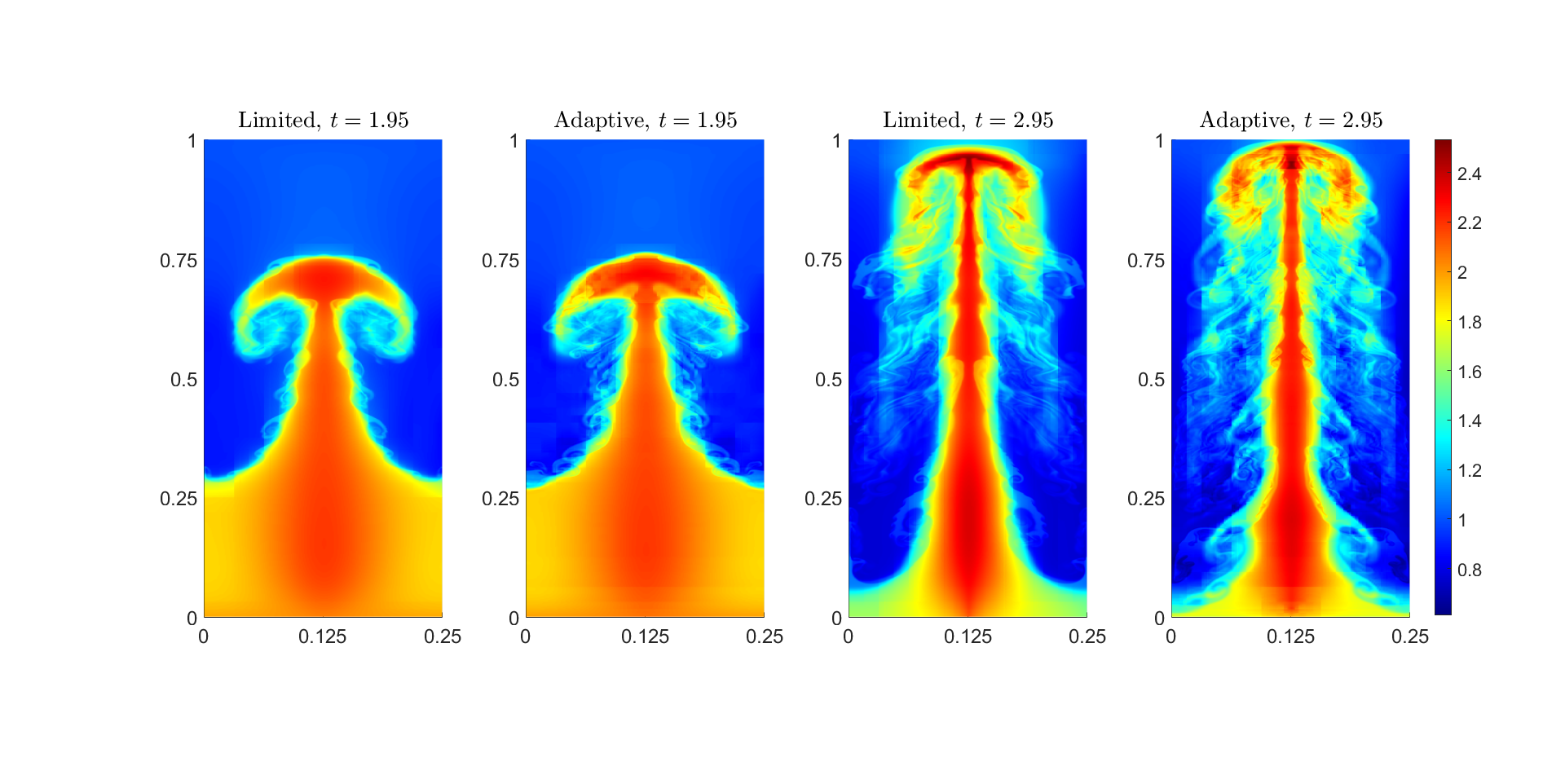}}
\caption{\sf Example 9:  Ces\`aro averages of the density $\rho^{\rm C}(1/2^{10})$ computed by the limited and adaptive schemes at $t=1.95$
and 2.95.\label{fig14bb}}
\end{figure}

As in Example 8, we also approximate the dissipative weak solution using the Ces\`aro averages computed by \eref{5.4} with the same sequence
of meshes. We present $\rho^{\rm C}(1/2^{10})$ computed by the limited and adaptive schemes in Figure \ref{fig14bb} at the times $t=1.95$
and 2.95. Once again, one can observe that the adaptive scheme better resolves the limiting dissipative weak solution.
\begin{rmk}
In this example, the solution is symmetric with respect to the vertical axis $x=0.125$. In order to enforce this symmetry, we have applied
the strategy from \cite{WDGK2020}: upon completion of each time evolution step, we replace the computed cell averages $\bm U_{j,k}$ with
$\widehat{\bm U}_{j,k}$, where
\begin{equation*}
\begin{aligned}
\widehat\rho_{j,k}&=\frac{\rho_{j,k}+\rho_{M-j,k}}{2},&(\widehat{\rho u})_{j,k}&=\frac{(\rho u)_{j,k}-(\rho u)_{M-j,k}}{2},\\
(\widehat{\rho v})_{j,k}&=\frac{(\rho v)_{j,k}+(\rho v)_{M-j,k}}{2},&\widehat E_{j,k}&=\frac{E_{j,k}+E_{M-j,k}}{2},
\end{aligned}
\end{equation*}
for all $j,k$ under the assumption that $j=1,\ldots,M$. Alternative symmetry enforcement techniques can be found in, e.g.,
\cite{DLGW,DLWG,Fleischman19,WTX}.
\end{rmk}

\begin{rmk}
As in the 1-D case, we also compare the computational costs of the studied limited and adaptive A-WENO schemes and present the CPU times
consumed by the adaptive scheme relative to the fully limited one. The obtained results are presented in Table \ref{tab3}, where one can see
that in the 2-D case, the difference in CPU times is slightly smaller than in the 1-D examples, but the adaptive scheme is still clearly more
efficient.
\begin{table}[ht!]
\centering
\begin{tabular}{|c|c|c|c|c|}
\hline
Example 5&Example 6&Example 7&Example 8&Example 9\\
\hline
80\% &81\% &81\% &80\% &83\%\\
\hline
\end{tabular}
\caption{\sf Example 5--9: CPU times consumed by the adaptive A-WENO scheme relative to the fully limited A-WENO scheme.\label{tab3}}
\end{table}
\end{rmk}

\section{Conclusion}\label{sec6}
In this paper, we have developed new adaptive alternative weighted essentially non-oscillatory (A-WENO) schemes for one- and two-
dimensional hyperbolic systems of conservation laws. The proposed schemes employ the scheme adaption strategy, according to which the
limited WENO-Z interpolation is only used to capture ``rough'' parts of the computed solution, while in the smooth areas, nonlimited fifth-
order interpolant is implemented. The ``rough'' regions are detected using a smoothness indicator. We have proposed a new, simple and
robust local smoothness indicator (LSI), which is based on the solutions computed at each of the three stages of the three-stage third-
order strong stability preserving Runge-Kutta time integrator. We have applied the new one- and two-dimensional adaptive A-WENO schemes to
the Euler equations of gas and dynamics using the recently proposed local characteristic decomposition (LCD) based central-upwind
numerical fluxes. We have conducted several numerical experiments and demonstrated that the new adaptive schemes are essentially non-
oscillatory and robust and, at the same time, more accurate than their fully limited counterparts.

In order to illustrate the high efficiency of the proposed adaptive A-WENO schemes, we have compared the CPU times consumed by the studied
fully limited and adaptive schemes in each numerical example. From the reported results, we conclude that the introduced scheme
adaption strategy leads to more efficient and, at the same time, more accurate A-WENO schemes. It should also be noted that if the
LCD-based CU numerical fluxes implemented in \eref{2.2f} and \eref{2.7f} are replaced with any other finite-volume (FV) numerical fluxes,
the resulting adaptive A-WENO schemes will still be substantially more efficient than the corresponding fully limited A-WENO schemes.
However, the difference in the CPU times may vary depending on the computational cost of the particular FV fluxes used.

\section*{Acknowledgments}
The work of A. Chertock was supported in part by NSF grants DMS-1818684 and DMS-2208438. The work of A. Kurganov was supported in part by
NSFC grants 12111530004 and 12171226, and by the fund of the Guangdong Provincial Key Laboratory of Computational Science and Material
Design (No. 2019B030301001).

\appendix
\section{The 1-D Fifth-Order WENO-Z Interpolant}\label{appa}
Here, we briefly describe the fifth-order WENO-Z interpolant.

Assume that the point values $W_j$ of a certain function $W(x)$ at the uniform grid points $x=x_j$ are available. We now show how to obtain
an interpolated left-sided value of $W$ at $x=x_\jph$, denoted by $W^-_\jph$. The right-sided value $W^+_\jph$ can then be obtained in the
mirror-symmetric way.

$W^-_\jph$ is computed using a weighted average of the three parabolic interpolants ${\cal P}_0(x)$, ${\cal P}_1(x)$ and ${\cal P}_2(x)$
obtained using the stencils $[x_{j-2},x_{j-1},x_j]$, $[x_{j-1},x_j,x_{j+1}]$, and $[x_j,x_{j+1},x_{j+2}]$, respectively:
\begin{equation}
W^-_\jph=\sum_{k=0}^2\omega_k{\cal P}_k(x_\jph),
\label{A1}
\end{equation}
where
\begin{equation}
\begin{aligned}
&{\cal P}_0(x_\jph)=\frac{3}{8}W_{j-2}-\frac{5}{4}W_{j-1}+\frac{15}{8}W_j,\\
&{\cal P}_1(x_\jph)=-\frac{1}{8}W_{j-1}+\frac{3}{4}W_j+\frac{3}{8}W_{j+1},\\
&{\cal P}_2(x_\jph)=\frac{3}{8}W_j+\frac{3}{4}W_{j+1}-\frac{1}{8}W_{j+2}.
\end{aligned}
\label{A2}
\end{equation}
Using a straightforward Taylor expansion one can show that \eref{A1}--\eref{A2} is fifth-order accurate if one takes the weights $\omega_k$
in \eref{A1} to be
\begin{equation}
\omega_k=\breve\omega_k:=\frac{d_k}{d_0+d_1+d_2},\quad d_0=\frac{1}{16},~~d_1=\frac{5}{8},~~d_2=\frac{5}{16},
\label{A3}
\end{equation}
resulting in the nonlimited point values, which we denote by
\begin{equation*}
\breve W^-_\jph:=\sum_{k=0}^2\breve\omega_k{\cal P}_k(x_\jph).
\end{equation*}

The computed interpolation may, however, be oscillatory in ``rough'' areas of $W(x)$ and thus the values $\breve W^-_\jph$ need to be
modified by replacing the weights \eref{A3} there with
\begin{equation}
\omega_k=\widetilde\omega_k:=\frac{\alpha_k}{\alpha_0+\alpha_1+\alpha_2},\quad
\alpha_k=d_k\left[1+\left(\frac{\tau_5}{\beta_k+\varepsilon}\right)^p\right],\quad\tau_5=|\beta_2-\beta_0|,
\label{A4}
\end{equation}
and thus obtaining the limited WENO-Z point values
\begin{equation}
\widetilde W^-_\jph:=\sum_{k=0}^2\widetilde\omega_k{\cal P}_k(x_\jph).
\label{A5}
\end{equation}
In \eref{A4}, $\beta_k$ are the following smoothness indicators for the corresponding parabolic interpolants ${\cal P}_k$:
\begin{equation}
\beta_k=\sum_{\ell=1}^2(\dx)^{2\ell-1}\int\limits_{C_j}\left(\frac{\partial^\ell{\cal P}_k}{\partial x^\ell}\right)^2{\rm d}x,\quad k=0,1,2.
\label{A6}
\end{equation}
Evaluating the integrals in \eref{A6}, we obtain
\begin{equation}
\begin{aligned}
&\beta_0=\frac{13}{12}\big(W_{j-2}-2W_{j-1}+W_j\big)^2+\frac{1}{4}\big(W_{j-2}-4W_{j-1}+3W_j\big)^2,\\
&\beta_1=\frac{13}{12}\big(W_{j-1}-2W_j+W_{j+1}\big)^2+\frac{1}{4}\big(W_{j-1}-W_{j+1}\big)^2,\\
&\beta_2=\frac{13}{12}\big(W_j-2W_{j+1}+W_{j+2}\big)^2+\frac{1}{4}\big(3W_{j}-4W_{j+1}+W_{j+2}\big)^2.
\end{aligned}
\label{A7}
\end{equation}
Finally, in all of the numerical examples reported in this paper, we have used $p=2$ and $\varepsilon=10^{-12}$.

\section{1-D Local Characteristic Decomposition}\label{appb}
Even though the WENO-Z interpolant \eref{A1}, \eref{A4}, \eref{A5}, \eref{A7} is essentially non-oscillatory, it is well-known that its
application to the conservative variables $\mU$ in a componentwise manner may lead to spurious oscillations in the computed solution. We, therefore, implement the reconstruction procedure described in Appendix \ref{appa} in the LCD framework.

Specifically, we first introduce the matrix $\widehat A_\jph:=A(\widehat\mU_\jph)$, where $\widehat\mU_\jph$ is either a simple average
$(\mU_j+\mU_{j+1})/2$ or another type of average of the $\mU_j$ and $\mU_{j+1}$ states (in the numerical examples reported in \S\ref{sec3}
and \S\ref{sec4}, we have used the simple average). As long as the system \eref{1.1} is strictly hyperbolic, we compute the matrices
$R_\jph$ and $R^{-1}_\jph$ such that $R^{-1}_\jph\widehat A_\jph R_\jph$ is a diagonal matrix and introduce the local characteristic
variables
in the neighborhood of $x=x_\jph$:
$$
\bm\Gamma_m=R^{-1}_\jph\mU_m,\quad m=j-2,\ldots,j+3.
$$
Equipped with the values $\bm\Gamma_{j-2}$, $\bm\Gamma_{j-1}$, $\bm\Gamma_j$, $\bm\Gamma_{j+1}$, $\bm\Gamma_{j+2}$, and $\bm\Gamma_{j+3}$,
we apply the interpolation procedure described in Appendix \ref{appa} to each of the components $\Gamma^{(i)}$, $i=1,\ldots,d$ of
$\bm\Gamma$ and obtain either the nonlimited $\breve{\bm\Gamma}_\jph^-$ or limited $\widetilde{\bm\Gamma}_\jph^-$ point values. The values
$\breve{\bm\Gamma}_\jph^+$ and $\widetilde{\bm\Gamma}_\jph^+$ are computed, as mentioned in Appendix \ref{appa}, in the mirror-symmetric
way. Finally, the corresponding nonlimited and limited point values of $\mU$ are given by
\begin{equation*}
\breve\mU^\pm_\jph=R_\jph\breve{\bm\Gamma}^\pm_\jph\quad\mbox{and}\quad\widetilde\mU^\pm_\jph=R_\jph\widetilde{\bm\Gamma}^\pm_\jph,
\end{equation*}
respectively.
\begin{remark}
A detailed explanation of how the average matrix $\widehat A_\jph$ and the corresponding matrices $R_\jph$ and $R^{-1}_\jph$ are computed in
the case of the Euler equation of gas dynamics can be found in, e.g., \cite{CCHKL_22}.
\end{remark}

\bibliographystyle{siam}
\bibliography{Chertock_Chu_Kurganov}
\end{document}